\documentclass[10pt]{amsart}

\usepackage{mathrsfs,dsfont,a4wide}
\usepackage{latexsym,amssymb,amsmath}
\usepackage{amsfonts}
\usepackage{eufrak}
\usepackage{graphicx}

\newcommand{\real}{\mathbb{R}}

\newcommand{\eps}{\varepsilon}
\newcommand{\cv}{\cos\varphi}
\newcommand{\sv}{\sin\varphi}
\newcommand{\ov}{\overline}
\newcommand{\vp}{\varphi}
\newcommand{\ad}{{\rm ad}}

\newcommand{\intt}{\int\limits}

\newcommand{\ham}{\vec h} 
 
\newcommand{\dd}{\partial}

\newcommand{\ppotR}[3]
{
\begin{figure}\begin{center}
~\includegraphics[height=#3truecm]{#1.eps}\\
\caption{#2}
\label{#1}
\end{center}
\end{figure}
\noindent$\!$}

\newcommand{\ppotRR}[5]
{

\begin{figure}\begin{center}
a)~\includegraphics[width=#4truecm,height=#5truecm]{#1.eps}
b)~\includegraphics[width=#4truecm,height=#5truecm]{#2.eps}\\
\caption{#3}
\label{#1}
\end{center}
\end{figure}
\noindent$\!$}

\newtheorem{theorem}{Theorem}
\newtheorem{Ex}{Example}

\newtheorem{defi}{Definition}       

\newtheorem{proposition}[theorem]{Proposition}
\newcommand{\BE}{\begin {equation}}
\newcommand{\EE}{\end {equation}}

\newcommand{\la}{\langle}
\newcommand{\ra}{\rangle}
\newtheorem{Theorem}{Theorem}

\begin{document}

\title[Sub-Riemannian minimal surfaces and structure of their singular sets]
{Minimal surfaces in sub-Riemannian manifolds\\and structure of their singular sets in the $(2,3)$ case.}
\author{Nataliya Shcherbakova}\address{SISSA/ISAS, via Beirut 2-4, 34100, Trieste, Italy, chtch$@$sissa.it} 

\date{}
\maketitle
\begin{abstract}

We study minimal surfaces in generic sub-Riemannian manifolds with sub-Riemannian structures of co-rank one. These surfaces can be defined  as the critical points of the  so-called {\it horizontal}  area functional associated to the canonical {\it horizontal} area form. We derive the intrinsic equation in the general case and then consider in greater detail $2$-dimensional surfaces in contact manifolds of dimension $3$. We show that in this case  minimal surfaces are projections of a special class of $2$-dimensional surfaces in the horizontal spherical bundle over the base manifold. Generic singularities of minimal surfaces turn out the singularities of this projection, and we give a complete local classification of them. We illustrate our results by examples in the Heisenberg group and the group of roto-translations. 

\end{abstract}

\section*{Introduction}

In the classical Riemannian geometry minimal surfaces realize the critical points of the area functional 
with respect to variations  preserving the boundary of a given domain. In this paper we study the generalization of the notion of minimal surfaces in sub-Riemannian manifolds known also as the Carnot-Carath\'eodory  spaces. 
This problem  was first introduced in the framework of Geometric Measure Theory for the Lie groups.
Mainly the obtained results (\cite{GarN}, ~\cite{GarP},~\cite{FSSC},~\cite{Pau0},~\cite{Malc},~\cite{Hw},~\cite{Rit}) concerns the Heisenberg groups, in particular ${\mathbb H}^1$; in  \cite{Citti} and \cite{Pau1} the authors were studying  the group $E_2$ of roto-translations of the plane, in \cite{Malc} there were also obtained some results for the case of  $S^3$. In \cite{Malc}, followed by just appeared paper  ~\cite{Yang},  the authors considered the problem in a more general setting and introduced the notion of minimal surfaces associated to CR structures in pseudohermitian manifolds of any dimension.

In this paper we develop a different approach using the methods of sub-Riemannian geometry.  
Though in particular cases of Lie Groups  ${\mathbb H}^m$, $E_2$ and $S^3$ the 
surfaces introduced in \cite{Malc} are minimal also in the sub-Riemannian sense, in general it is not true.
The sub-Riemannian point of view on the problem is based on the following construction.

Consider an $n$-dimensional smooth  manifold $M$ and a co-rank $1$ smooth vector distribution $\Delta$ in it (``horizontal'' distribution). It is assumed that the sections of $\Delta$ are endowed with a Euclidean structure, which can be described by fixing an orthonormal  basis of vector fields $X_1,\dots,X_{n-1}$ on $\Delta$ (see ~\cite{Bel}). Then $\Delta$ defines a {\it sub-Riemannian structure} in $M$. In this case $M$ is said to be a sub-Riemannian manifold. Given a sub-Riemannian structure there is a canonical way to define a volume form $\mu\in\Lambda^n M$ associated to it. In addition, for any hypersurface $W\subset M$ the horizontal unite vector $\nu$ such that  
$$
\intt_\Omega i_{\nu}\mu=\max_{{X\in \Delta}\atop{\|X\|_{\Delta}=1}}\intt_\Omega i_X\mu,\qquad \Omega\subset W,
$$
plays the role of the Riemannian normal in the classical case, and the $n-1$-form $i_{\nu}\mu$  defines  the {\it horizontal area form} on $W$.  All these notions are  direct generalizations of the the classical ones in the Riemannian geometry (i.e., in the case $\Delta\equiv TM$).  
 
Going further in this direction, we  define sub-Riemannian minimal surfaces in $M$  as the critical points
of the  functional associated to the horizontal area form.  It turns out that these surfaces 
satisfy  the following  intrinsic equation
\BE\label{eq0}
(d\circ i_{\nu}\mu)\Big|_{W\setminus\Sigma}=0,
\EE
where $\Sigma$ is the set of the so-called {\it characteristic points} of $W$, i.e., the  points where $W$ is tangent to $\Delta$.  The described construction does not require the existence of any additional global structure in $M$, and can be generalized for sub-Riemannian structures of greater co-rank. 
 
\vspace{5 pt}

The existence of the singular set $\Sigma$ is one of the main difficulties of the problem. 
In general, the set $\Sigma$ can be quite large and have its own non-trivial intrinsic geometry.
In the second part of this paper we show how this problem can be resolved in the case of $2$-dimensional surfaces in  $3$-dimensional contact manifolds.

It turns out that in the $(2,3)$ case, due to the relatively small dimension, there is an elegant way to extend the definition of a sub-Riemannian minimal surface over its singular set. Namely, in this case all information related to the intrinsic geometry of a surface $W$ is encoded in its {\it characteristic curves} $\gamma:\,[0,T]\mapsto W$ such that  $\dot\gamma(t)\in T_{\gamma(t)}W\cap \Delta_{\gamma(t)}$ for all $t\in [0,T]$.  The vector field $\eta$ 
(the {\it characteristic vector field}) tangent to characteristic curves is $\Delta$-orthogonal to the sub-Riemannian normal of $W$ and  it is well defined (as well as $\nu$ and the sub-Riemannian area form) away from the set $\Sigma$, where 
$\nu$ degenerates. On the other hand, $\eta$ turns to be a projection onto $M$ of a special invariant vector filed $V$  in the horizontal spherical bundle $S_{\Delta}M$ over $M$. In contrast with $\nu$, the vector filed $V$ is well defined everywhere in $S_{\Delta}M$, and moreover, minimal surface equation (\ref{eq0}) can be transformed into a quasilinear equation whose characteristics are exactly the integral curves of $V$. 

These observations motivate the key idea of our work. Namely, we consider the sub-Riemannian minimal surfaces in the $(2,3)$ case  as the projections of  $2$-dimensional surfaces foliated by the integral curves of $V$ (the {\it generating surfaces} of the sub-Riemannian minimal surfaces). So instead of dealing with a highly degenerate PDE (\ref{eq0}), we
study solutions of the ODE associated to the vector field $V$: $\dot{\bar q}=V(\bar q)$, 
$\bar q\in \ov M\simeq S_{\Delta}M$. By varying initial conditions, we provide a local characterization of all possible sub-Riemannian minimal surfaces, together with their singular sets. 

From the point of view that we develop in this paper characteristic points of sub-Riemannian minimal surfaces are either singular points of their generating surfaces or the singularities of the projection of the last ones onto the 
base manifold $M$. Thus among characteristic points of the minimal surfaces generated by regular surfaces 
in $S_{\Delta}M$ we distinguish between {\it regular} and {\it singular} characteristic points. 
The first ones are regular points of the projection of the generating surface,  and they form {\it simple singular curves}. The generic situation for singular characteristic points, i.e., the singular point of the projection, is described in Theorem 1:

\vspace{5pt}
\noindent{\bf Theorem 1} {\it Let $\bar q$ be a regular point of a  generating surface ${\mathcal W}$ such that $q=\pi[\bar q]$ is a  generic singular  characteristic point of $W=\pi[{\mathcal W}]$. Then:\\
a) a small enough  neighborhood of $q$ in $W$ contains a pair of simple singular curves, which touch each other at $q$ and form a unique smooth curve on $W$ passing through $q$;\\
b) there exists a choice of local coordinates $(\tau,\sigma)\in\Omega\subseteq \real^2$ such that   
$$
W=\{q(\tau,\sigma)=(\tau^2,\sigma,\tau\sigma),\;\; \tau, \sigma\in \Omega\subseteq\real^2\},\qquad q=q(0,0), 
$$
i.e., in the neighborhood of $q$ $W$ has the structure of Whitney's umbrella.}
 \vspace{5pt}

In less generic situations other types of singularities may appear. For instance, certain minimal surfaces may contain isolated characteristic points or entire curves of singular points ({\it strongly singular curves}). All described types of singularities are already present in the Heisenberg group ${\mathbb H}^1$,  which we use to illustrate our results.

The author is grateful to prof. A. Agrachev whose vision of the problem inspired  this work.

\section{Minimal surfaces in sub-Riemannian manifolds}

\subsection{Sub-Riemannian structures and associated objects}   

Let $M$ be an $n$-dimensional smooth manifold. Consider a co-rank $1$ vector distribution  $\Delta$ in $M$: 
$$
\Delta=\bigcup_{q\in M}\Delta_q,\quad \Delta_q\subset T_qM,\quad q\in M.
$$
By definition, the sub-Riemannian structure in $M$  is a pair $(\Delta, \langle\cdot,\cdot\rangle_{\Delta})$, where
$\langle\cdot,\cdot\rangle_{\Delta}$ denotes a  smooth family of Euclidean inner products on $\Delta$.
In what follows we will call $\Delta$ {\it the horizontal distribution} and  keep 
the same notation $\Delta$ both for the vector distribution and for the associated sub-Riemannian structure.

Let $X_i$, $i=1,\dots, n-1$, be a horizontal orthonormal basis:
$$
\Delta_q={\rm span}\{X_1(q),\dots, X_{n-1}(q)\},\qquad q\in M,
$$
$$
\la X_i(q),X_j(q)\ra_{\Delta}=\delta_{ij},\qquad q\in M,\;i,j=1,\dots, n-1.
$$
By $\Theta\in \Lambda^{n-1}\Delta$ we will denote the Euclidean volume form on $\Delta$.

In what follows we will assume that $\Delta$ is {\it bracket-generating}. In the present case this means that
$$
{\rm span}\{X_i(q), [X_i, X_j](q),\;i,j=1,\dots,n-1,\;q\in M\}=T_q M.
$$
Hereafter the square brackets denote the Lie brackets of vector fields.
If $\Delta$ is bracket-generating, then by the Frobenius theorem it is completely non-holonomic, i.e.,
there is no invariant sub-manifold in $M$ whose tangent space coincides with 
$\Delta$ at any point.

We can also define the  distribution $\Delta$ as the kernel of some differential $1$-form. 
Let  $\omega\in \Lambda^1M$ be such a form : 
$$
\Delta_q={\rm Ker}\, \omega_q=\{v\in T_q M:\;\omega_q(v)=0\},\quad q\in M.
$$
It  is easy to check that $\Delta$ is bracket-generating at $q\in M$ 
if and only if $d_q\omega\ne 0$.

Though in general the form $\omega$  is defined up to a multiplication by a non-zero scalar function, there is a canonical way to choose it by using the Euclidean structure on $\Delta$.  Indeed, by standard construction the Riemannian structure on $\Delta$ can be extended to the spaces of forms $\Lambda^k \Delta$, $k\le n-1$.
In particular, for any  $2$-form $\sigma$ we set
$$
\|\sigma_q\|_{\Delta}=\left(\sum_{i,j=1\atop i<j}^{n-1} \sigma_q(X_i(q), X_j(q))^2\right)^{\frac{1}{2}},
$$
$\{X_i(q)\}_{i=1}^{n-1}$, as before, being  the orthonormal horizontal basis of $\Delta_q$.
Now we  can normalize $\omega$ as follows:
\BE\label{omega}
\omega_q(\Delta_q)=0,\qquad \|d_q\omega\|^2_{\Delta}=1,\qquad \forall q\in M.
\EE
The $1$-form satisfying (\ref{omega}) is the  {\it canonical $1$-form associated to $\Delta$}.
In the fixed horizontal orthonormal basis $\{X_i(q)\}_{i=1}^{n-1}\in \Delta_q$ equations (\ref{omega}) become
\BE\label{omega1}
\omega_q(X_i(q))=0,\qquad \sum_{i,j=1\atop i<j}^{n-1} d_q\omega(X_i(q), X_j(q))^2=1,\qquad i=1,\dots, n-1.
\EE
Clearly the canonical $1$-form $\omega$ is defined up to a sign and does not depend on the choice of the horizontal basis.
In local coordinates on $M$ its components can be expressed in terms of  the coordinates of the vector 
fields $X_i$ and their first derivatives. Indeed, we have
$$
\|d_q\omega\|_\Delta^2=\sum_{i,j=1\atop i<j}^{n-1} d_q\omega(X_i(q), X_j(q))^2=\sum_{i,j=1\atop i<j}^{n-1}\omega_q([X_i,X_j](q))^2
$$
because according to Cartan's formula for any pair of vector fields $X$ and $Y$
$$
d\omega(X,Y)=X \omega(Y)-Y \omega(X)-\omega([X,Y]).
$$

Once the orientation in $M$ if fixed by a choice of the sing of $\omega$, the volume form
$$
\mu=\Theta\wedge \omega
$$
is uniquely defined. We will call this volume form {\it the canonical volume form associated to $\Delta$}. 
The canonical volume form  $\mu$ is a ``global'' object  in $M$, though it is intrinsically defined by the sub-Riemannian structure $\Delta$.

\subsection{Horizontal area form}

Let $W\subset M$, ${\rm dim} W=n-1$ be a smooth hyper-surface in $M$ and let $\Omega\subset W$ be an open domain.
For simplicity we assume that the  vector field $X\in TM$ is transversal to $W$, though this assumption is not restrictive.
Let us denote by $e^{t X}$ the flow  generated by $X$ in $M$ and consider the map
$$
\Pi^X:\qquad [0,\varepsilon]\times \Omega\mapsto M,
$$
$$
\Pi^X(t,q)=e^{t\,X}(q),\quad q\in M.
$$
Denote by 
\BE
\Pi_{(\eps, \Omega)}^X=\left\{e^{tX}(q),\;q\in\Omega,\;t\in[0,\varepsilon]\right\} 
\EE
the cylinder formed by the images of $\Omega$ translated along the integral curves of $X$ 
parametrized by $t\in[0,\eps]$. Clearly, $\Pi_{(0, \Omega)}^X\equiv\Omega$.
By definition,
$$
Vol(\Pi_{(\eps, \Omega)}^X)=\int\limits_{\Pi^X_{(\eps, \Omega)}}\mu =\int\limits_{[0,\varepsilon]\times\Omega}({\Pi^X})^*\mu,
$$  
where $({\Pi^X})^*$ denotes the pull-back map associated to $\Pi^X$ and $\mu$ is the canonical volume form defined above\footnote{Here  
we use the canonical volume form associated to $\Delta$, though the whole construction works for any volume form in $M$.}.

\begin{defi} We will call 
\BE\label{areaform}
A_\Delta(\Omega)=\max_{{X\in \Delta}\atop{\|X\|_{\Delta}=1}}\;\lim_{\varepsilon\to 0}\frac{Vol(\Pi^X_{(\eps,\Omega)})}{\eps}
\EE
the sub-Riemannian  (or {\it horizontal}) area of the domain $\Omega$ associated to $\Delta$.\end{defi}

\noindent{\bf Remark.} The horizontal area (\ref{areaform}) is nothing but the generalization of the classical notion of the Euclidean area:  it defines the area of the base of a cylinder as the ratio of its volume and height.  

\vspace{5pt}

Let us find a more convenient expression for (\ref{areaform}). First of all observe that since 
$({\Pi^X})^*\mu$ is a form of maximal rank $n$ in $M$ it follows that $dt\wedge({\Pi^X})^*\mu=0$. Hence
$$
0=i_{\dd_t}\left(dt\wedge({\Pi^X})^*\mu\right)=i_{\dd_t}dt\wedge({\Pi^X})^*\mu-dt\wedge i_{\dd_t}({\Pi^X})^*\mu,
$$
i.e.,
$$
({\Pi^X})^*\mu=dt\wedge i_{\dd_t}({\Pi^X})^*\mu.
$$
Taking into account that $\Pi^X_*\dd_t=X$ we obtain
$$
A_\Delta(\Omega)=\max_{{X\in \Delta}\atop{\|X\|_{\Delta}=1}}\;\lim_{\varepsilon\to 0}\frac{Vol(\Pi^X_{(\eps,\Omega)})}{\eps}=
$$
$$
= \max_{{X\in \Delta}\atop{\|X\|_{\Delta}=1}}\;\frac{\dd}{\dd \eps}\Big|_{\eps=0}\intt_{\Pi^X_{(\eps,\Omega)}}\mu
=\max_{{X\in \Delta}\atop{\|X\|_{\Delta}=1}}\;\frac{\dd}{\dd \eps}\Big|_{\eps=0} \intt_{[0,\eps]\times\Omega}dt\wedge i_{\dd_t}({\Pi^X})^*\mu=
$$
$$
= \max_{{X\in \Delta}\atop{\|X\|_{\Delta}=1}}\;\frac{\dd}{\dd \eps}\Big|_{\eps=0}\intt_0^{\eps}\Big(\intt_{\Pi^X_{(t,\Omega)}}i_X\mu\Big) dt=
 \max_{{X\in \Delta}\atop{\|X\|_{\Delta}=1}}\intt_\Omega i_X\mu.
$$

\begin{defi}\label{normal}
The horizontal unite  vector field $\nu\in\Delta$, $\|\nu\|_{\Delta}=1$, such that
$$
\intt_\Omega i_{\nu}\mu=\max_{{X\in \Delta}\atop{\|X\|_{\Delta}=1}}\intt_\Omega i_X\mu
$$
is called the sub-Riemannian  or horizontal normal of $\Omega\subset W$. 
The $(n-1)$-form $i_{\nu}\mu$ is called the sub-Riemannian  or horizontal area form on $W$ associated to $\Delta$.
\end{defi}

\vspace{5 pt}

\noindent{\bf Remark} According to the definition, the horizontal normal $\nu$ is well defined
everywhere except the points where the surface $W$ it is tangent to the distribution $\Delta$. Such points are called the {\it characteristic points} of $W$ and they form the subspace 
$$
\Sigma=\{q\in W:\quad T_q W=\Delta_q\}
$$
called the {\it singular set} of $W$. The set $\Sigma$ can have a very non-trivial intrinsic geometry.
In Section 2  we will  analyze in detail the structure of $\Sigma$ in the case ${\rm dim}M=3$.

\vspace{5 pt}

The sub-Riemannian normal is an intrinsic object associated to any hypersurface in $M$, and Definition \ref{normal} does not require any global structure in $M$. Nevertheless,  if $M$ is a Riemannian manifold whose Riemannian structure is compatible with the sub-Riemannian structure on $\Delta$, i.e., if the inner product $\la\cdot, \cdot\ra$ on $TM$ satisfies $\la\cdot,\cdot\ra_{\Delta}=\la\cdot, \cdot\ra\big|_\Delta$, then it is easy to see  that the sub-Riemannian normal $\nu$ is  nothing but the projection on $\Delta$ of the Riemannian  unit normal $N$ of $W$, normalized w.r.t. $\|\cdot\|_{\Delta}$. This fact follows from the relation
$$
\int\limits_{\Omega} i_{X} \mu =\int\limits_{\Omega}\langle X,N\rangle i_N\,\mu,\qquad \forall  X\in Vec(M).
$$
Thus if $X_1,\dots,X_{n-1}\in \Delta$ is an orthonormal horizontal basis of $\Delta$, then 
\BE\label{srnorm}
\nu=\sum\limits_{i=1}^{n-1}\nu_i X_i\,,\qquad \nu_i=\frac{\langle N,X_i\rangle}
{\sqrt{\langle N,X_1\rangle ^2+\dots+\langle N,X_{n-1}\rangle ^2}},
\EE
and the  horizontal area form reads
$$
i_{\nu}\mu=\langle \nu, N \rangle i_N \,\mu=
\sqrt{\langle N,X_1\rangle ^2+\dots+\langle N,X_{n-1}\rangle ^2}\,i_N\,\mu.
$$
Consider now a hyper-surface $W$ defined as a level set of a smooth, let us say $C^2$, function:
$$
W=\left\{q\in M:\quad F(q)=const,\;F\in C^2(M), \;d_q F\ne 0\right\}.
$$
If $X \equiv X_n$ is a vector field transversal to $W$ and such that $\{X_i(q)\}_{i=1}^n$ 
form an orthonormal basis of $T_q M$ at $q\in M$, then 
$$
N(q)=D_0^{-1}\sum\limits_{i=1}^n X_i F(q)\,X_i(q),\qquad D_0=\left(\sum\limits_{i=1}^n X_i F(q)^2\right)^{1/2}
$$
and
\BE\label{nu}
\nu(q)=D_1^{-1}\sum\limits_{i=1}^{n-1} X_i F(q)\, X_i(q) ,\qquad D_1=\left(\sum\limits_{i=1}^{n-1} X_i F(q)^2\right)^{1/2}.
\EE
Hereafter  $X_iF$ denotes the directional derivative of the function  $F$ along the vector field $X_i$.

\subsection{Sub-Riemannian minimal surfaces}

Let us  compute the first  variation of the horizontal area $A_\Delta(\cdot)$.
Take a bounded domain $\Omega\subset W$ and a vector filed $V\in Vec(M)$ such that $V\big|_{{\partial}\Omega}=0$.
For the moment we assume that $\Omega$ contains no characteristic points.
Consider a one-parametric family of hyper-surfaces generated by the vector field $V$
$$
\Omega^t= e^{t V}\Omega,\qquad \Omega^0\equiv\Omega,
$$
and denote by $\nu^t$ the horizontal unit normals to $\Omega^t$.
We have
$$
A_{\Delta}(\Omega^t)=\int\limits_{e^{t V}\Omega} i_{\nu^t}\ \mu=\int\limits_{\Omega}(e^{t V})^* i_{\nu^t}\,\mu
=\int\limits_{\Omega}e^{t L_V} i_{\nu^t}\,\mu.
$$
Further, 
\BE\label{int1}
\frac{\partial}{\partial t} \Big|_{t=0} A_{\Delta}(\Omega^t)=\int\limits_{\Omega}L_V i_{\nu}\,\mu+
\int\limits_{\Omega}i_{\frac{\partial \nu^t}{\partial t}\big|_{t=0}}\mu.
\EE
It is not hard to show that the second integral in (\ref{int1}) vanishes, 
because the horizontal vector field  $\frac{\dd \nu^t}{\dd t}\big|_{t=0}$ is tangent to $\Omega$.
Indeed, at any generic (non-characteristic) point $q\in \Omega$ we have $\nu(q)\notin T_q \Omega$ and 
${\rm dim}\Delta_q\cap T_q \Omega=n-2$. On the other hand,  differentiating 
the equality $\la \nu^t,\nu^t\ra_{\Delta}=1$ we get 
\BE\label{unity}
\langle \frac{\partial \nu^t}{\partial t}\Big|_{t=0},\nu\rangle_{\Delta}=0,
\EE
and hence  $\frac{\dd \nu^t}{\dd t}\big|_{t=0}(q)\in T_q \Omega$.

Using Cartan's formula we can transform the first part of (\ref{int1}) as follows:
$$
\int\limits_{\Omega}L_V i_{\nu}\,\mu=\int\limits_{\Omega}(i_V\circ d+d\circ i_V) i_{\nu}\,\mu=
\int\limits_{\Omega}(i_V\circ d \circ i_{\nu})\mu+\int\limits_{\Omega}(d\circ i_V\circ i_{\nu})\mu.
$$
Applying the Stokes theorem to the second integral we see that it vanishes:
$$
\int\limits_{\Omega}(d\circ i_V\circ i_{\nu})\mu=\int\limits_{\partial \Omega}(i_V\circ i_{\nu})\mu=0
$$
provided $V\big|_{\partial \Omega}=0$ and $\dd \Omega$ is sufficiently regular.
Thus, 
$$
\frac{\partial}{\partial t} \Big|_{t=0} A_{\Delta}(\Omega^t)=\int\limits_{\Omega}i_V\circ(d \circ i_{\nu}\mu).
$$
 
\begin{defi}\label{MinS} We  say that the hyper-surface $W$ is a minimal surface w.r.t. the sub-Riemannian structure $\Delta$ (or just $\Delta$-minimal) iff     
\BE\label{mineq}
( d\circ i_{\nu}\mu)\Big|_{W\setminus\Sigma}=0,
\EE
where $\Sigma$ is the singular set of $W$.
\end{defi}
  
Needless to say that the minimality property of a hypersurface does not depend on the chosen orientation in $M$.
The whole construction can be further generalized for the case of vector distributions of co-rank greater than $1$.

\subsection{Canonical form of the minimal surface equation in contact sub-Riemannian  manifolds}

The construction described above works in any sub-Riemannian manifold endowed with a co-rank $1$ distribution.
From now on we restrict ourselves to the case $n=2m+1$ and assume that the distribution $\Delta$ is contact, 
i.e., the $2m+1$-form $(d\omega)^m\wedge \omega$ is non-degenerate. In this case we say that $M$  is a {\it contact sub-Riemannian manifold}.

First of all we recall that in the contact case there exists a special  uniquely defined vector filed $X\in TM$ associated to $\omega$.  This vector field is called {\it the Reeb vector field} of the contact form $\omega$ and it satisfies the following equalities 
\BE\label{Reeb}
\omega_q(X(q))=1, \qquad d_q\omega(v, X(q))=0,\qquad \forall v\in \Delta_q.
\EE
Using this vector field we can canonically extend the sub-Riemannian structure on $\Delta$ to the whole $TM$. The resulting Riemannian structure in $M$ is compatible with $\Delta$ by construction. The basis of vector fields 
$\{X_1,\dots, X_{2m},X\}$ is then a canonical basis associated to the contact sub-Riemannian structure $\Delta$.

\vspace{5 pt}

Set  $X_{2m+1}\equiv X$ and denote by $c_{ij}^k\in C^\infty(M)$ {\it the  structural constants} of the canonical frame $\{X_i\}_{i=1}^{2m+1}$:
\BE\label{se}
[X_i, X_j]=-\sum\limits_{k=1}^{2m+1} c_{ij}^k X_k.
\EE
Let $\{\theta_i\}_{i=1}^{2m+1}$ be  the basis of $1$-forms dual to $\{X_i\}_{i=1}^{2m+1}$. 
Clearly, $\theta_{2m+1}\equiv\omega$ and the canonical volume form  reads
$$
\mu=\theta_1\wedge\dots\wedge\theta_{2m+1}.
$$
We also recall that from Cartan's formula it follows that
\BE\label{Cart}
d\theta_k=\sum_{i,j=1\atop i<j}^{2m+1}c_{ij}^k\theta_i\wedge\theta_j,\qquad k=1,\dots, 2m+1. 
\EE
 
Now we are ready to derive the canonical form of the minimal surface equation (\ref{mineq}) in contact case. We have 
$$
i_\nu\,\mu=\left(\sum_{k=1}^{2m}(-1)^{k+1}\nu_k \,\theta_1\wedge\dots\wedge\widehat{\theta_k}\wedge\dots\wedge \theta_{2m}\right)\wedge\theta_{m+1}=\Xi\wedge\theta_{2m+1}.
$$
Here $\widehat{\theta_k}$ denotes the omitted element in the wedge product and 
$$\Xi=\sum_{k=1}^{2m}(-1)^{k+1}\nu_k \,\theta_1\wedge\dots\wedge\widehat{\theta_k}\wedge\dots\wedge \theta_{2m}.$$
Further, 
$$
d\,i_\nu\,\mu=d \Xi\wedge\theta_{2m+1}-\Xi\wedge d\theta_{2m+1}.
$$
Recalling now that $d\nu_k=\sum\limits_{i=1}^{2m+1}\,X_i\nu_k\,\theta_i$, we obtain
$$
d\Xi\wedge\theta_{2m+1}=\sum_{k=1}^{2m}(-1)^{k+1}\left(d\nu_k\wedge \theta_1\wedge\dots\wedge\widehat{\theta_k}\wedge\dots\wedge \theta_{2m}+\right.
$$
$$
\left.+\nu_k\, d(\theta_1\wedge\dots\wedge\widehat{\theta_k}\wedge\dots\wedge \theta_{2m})\right)\wedge\theta_{2m+1}=
\left(\sum_{k=1}^{2m}X_k\nu_k+\sum_{j=1}^{2m}\nu_k c_{kj}^j\right)\mu.
$$
On the other hand, 
$$
\Xi\wedge d\theta_{2m+1}=\Xi\wedge \sum_{i,j=1\atop i<j}^{2m+1}c_{ij}^{2m+1}\theta_i\wedge\theta_j=-\left(\sum_{k=1}^{2m}\nu_k c_{k 2m+1}^{2m+1}\right)\mu.
$$
Summing up we obtain the following equation:
\BE\label{mineqN}
\left.\left[{\rm div}^{\Delta}\nu+\sum\limits_{i=1}^{2m} \nu_i\left(\sum\limits_{j=1}^{2m+1} c_{ij}^j\right)\right]
\right|_{W\setminus\Sigma}=0. 
\EE
The left-hand side of (\ref{mineqN}) corresponds to the {\it sub-Riemannian mean curvature} of the hyper-surface $W$, while its first term
$$
{\rm div}^{\Delta}\nu=\sum_{i=1}^{2m} X_i\nu_i
$$
is called {\it the horizontal divergence} of the sub-Riemannian normal $\nu$. 
Equation (\ref{mineqN}) is the canonical equation of sub-Riemannian minimal surfaces in a contact sub-Riemannian manifold.

\vspace{5pt}

\noindent{\bf Remark} In \cite{Malc} there was introduced the notion of minimal surfaces associated to CR structures in pseudohermitian manifolds. This approach later  was used in \cite{Hw}, \cite{Yang}, and in \cite{Pau1}.  In general the sub-Riemannian structures we consider  in the present paper are not equivalent to the CR structures, and  the class of surfaces satisfying  (\ref{mineqN}) differs from its analog defined in \cite{Malc}. Nevertheless, in many particular cases,  for instance, in the cases of sub-Riemannian structures associated to the Heisenberg group, the group of roto-translations,  and $S^3$, these structures coincide. Thus in these cases our results are comparable with the ones of cited papers.

\label{examp1}\begin{Ex}(The Heisenberg distribution) {\rm
Let $M=\real^{2m+1}$ and denote by  $(x_1,\dots,x_{2m},t)=q$ the Cartesian coordinates in $M$.
Let  $\Delta$ be such that $\Delta_q={\rm span}\{X_i(q)\}_{i=1}^{2m}$, $q\in M$, where
\BE\label{heis}
X_i(q)=\dd_{x_i}-\frac{x_{i+m}}{2}\dd_z,
\EE
$$
X_{i+m}(q)=\dd_{ x_{i+m}}+\frac{x_i}{2}\dd_t,\quad i=1,\dots,m.
$$
The vector distribution $\Delta$ is characterized by the following commutative relations:
\BE\label{deltaH}
[X_i,X_j]=0,\quad {\rm for}\quad j\ne i+m,\qquad [X_i,X_{i+m}]=\dd _t,
\EE
and therefore it is a co-rank $1$ bracket-generating distribution.
The vector fields $X_i$, $i=1,\dots,2n$, generate the {\it Heisenberg Lie algebra} in $\real^{2m+1}$.
In what follows we will call any vector  distributions, which satisfy the commutative relations (\ref{deltaH}), 
the {\it Heisenberg distribution}. The space $\real^{2m+1}$ endowed with the structure of this distribution is the
Heisenberg group ${\mathbb H}^{m}$.

From (\ref{omega1}) we find the canonical $1$-form $\omega$: 
\BE\label{omegaH}
\omega=\pm\frac{1}{\sqrt{m}}(dt -\frac{1}{2}\sum_{i=1}^m(x_{i}\,dx_{i+m}-x_{i+m}\,dx_{i})).
\EE
Clearly $\omega$ is a contact form since 
$(d\omega)^m\wedge d\omega=\pm\frac{1}{m^m}\bigwedge\limits_{i=1}^{2m}dx_i\wedge dt$ is non-degenerate. 
The associated Reeb vector field is $X=\pm \sqrt m \dd_t$.  The only non-zero structural constants of the canonical frame are   $c_{i i+m}^{2m+1}=\pm \frac{1}{\sqrt m}$, $i=1,\dots,m$.
Due to the high degeneracy of the sub-Riemannian structure the canonical 
minimal surface equation takes a very simple form:   
$$
{\rm div}^{\Delta}\nu\;\Big|_{W\setminus\Sigma}=0.
$$
This is the well known minimal surface equation in the Heisenberg group (see \cite{GarP}, \cite{FSSC}, \cite{Malc}, \cite{Rit}, etc.)
}\end{Ex}

\section{Sub-Riemannian minimal surfaces associated to the $(2,3)$ contact vector distributions}

The less dimensional situation where sub-Riemannian minimal surfaces appear is the case of a contact distribution $\Delta$ of rank $2$ in the $3$-dimensional manifold $M$. Almost all known results on sub-Riemannian minimal surfaces are related to this case  and concern minimal surfaces in Lie groups mentioned above. In this paper we try to give a general picture.
Our main idea is to study (\ref{mineqN}) using the classical method of characteristics. This point of view permits to describe all possible minimal surfaces,  and, as we will show in a while,  to resolve the problem of the presence of singular sets.  All this is possible due to the fact that  ${\rm dim}T_qM\cap\Delta_q=1$ at any non-characteristic point $q\in M$.

\subsection{The $(2,3)$ structures}
 
From now on we set $n=3$ and $\Delta=\bigcup_{q\in M}\Delta_q$ with $\Delta_q={\rm span}\{X_1(q)$, $X_2(q)\}$.
As before, we complete the horizontal frame by the Reeb vector field $X\equiv X_3$ associated to $\Delta$,
and denote by $c_{ij}^k$ the structural constant of the canonical frame $\{X_i\}_{i=1}^3$. By definition,
 $c_{ij}^k=-c_{ji}^k$. Moreover, (\ref{Reeb}) and (\ref{se}) imply
\BE\label{se3}
c_{12}^3=\pm1,\qquad c_{13}^3=c_{23}^3=0.
\EE
More symmetry  relations of the structural constants can be obtained from the 
Jacobi identity
$$
[X_1,[X_2,X_3]]+[X_3,[X_1,X_2]]+[X_2,[X_3,X_1]]=0.
$$
In particular, if $M$ is a Lie group, the structural constants do not depend on the points of the base manifold $M$, and
the Jacobi identity implies the additional symmetry relations:
\BE\label{jac}
c_{13}^1+c_{23}^2=0,\qquad
c_{12}^1c_{13}^1+c_{12}^2c_{23}^1=0,\qquad
c_{12}^1c_{13}^2+c_{12}^2c_{23}^2=0.
\EE

\vspace{5 pt}

Let us consider a regular hyper-surface $W\subset M$  and let $\nu\in \Delta$ be its horizontal normal.  
As we know, $\nu$ is well defined on $W\setminus\Sigma$, where $\Sigma$ is the singular set of $W$, which contains characteristic points of $W$. If $W$ is $\Delta$-minimal, then by (\ref{se3})  
\BE\label{mineq1}
\Big(X_1\nu_1+X_2\nu_2+\nu_1 c_{12}^2 - \nu_2 c_{12}^1\Big)\Big|_{W\setminus\Sigma}=0.
\EE
Assume that $W$ is  given as a level set of some smooth function $F$. Then $\Sigma=\{q\in M:\;X_1 F(q)=X_2 F(q)=0\}$.
By (\ref{nu}), away from $\Sigma$ the function $F$ should satisfy the following PDE:
\BE\label{mineqF}
\Big[\Big(X_1^2F\,(X_2 F)^2+X_2^2 F\,(X_1F)^2-X_1F\, X_2 F\,(X_1\circ X_2+X_2\circ X_1)F\Big)D_1^{-3}+
\EE
$$
+\Big(c_{12}^2 X_1F-c_{12}^1X_2F\Big)D_1^{-1}\Big]\Big|_{W\setminus \Sigma}=0,\qquad D_1=\sqrt{X_1F^2+X_2F^2}. 
$$
Some non-trivial solutions of this equation are known, especially in the particular case of ${\mathbb H}^1$, the interested reader can consult in \cite{Malc} and other papers from the Bibliography. Let us also consider another important for applications case of the distribution associated to the  Lie Group $E_2$ (the group of rotation and translations of the plane). 

\label{examp2}\begin{Ex}{\rm The Lie group $E_2$ can be realized as $\real^2\times {\mathbb S}^1$.
In local coordinates $(x,y,z)$ with $(x,y)\in \real^2$, $z\in{\mathbb S}^1$  the left-invariant basis of the  corresponding Lie algebra  is given by   
vector fields 
$$
X_1=\cos z \dd_x+\sin z\dd_y,\qquad  X_2=\dd_z.
$$
It is easy to check that the  horizontal distribution  $\Delta^{E^2}$ with sections 
$\Delta^{E^2}_q={\rm span}\{X_1(q),X_2(q)\}$, $q\in M$, is contact, 
the corresponding canonical $1$-form is $\omega=\pm(\sin z dx -\cos z dy)$. 
The Reeb vector field coincides with the Lie bracket $[X_1,X_2]$ (up to the sign) and the only non-zero 
structural constants are $c_{12}^3=c_{23}^1=\pm 1$. Therefore, the minimal surface equation, as in the Heisenberg case, contains only the  divergence  term:
$$
(X_1\,\nu_1+X_2 \nu_2)\Big|_{W\setminus\Sigma}=0.
$$
For instance one can easily check that the following surfaces are $\Delta^{E^2}$-minimal:\\
a). $y=x+ B(\sin z+\cos z)+C,\qquad B,C={\rm const}$;\\
b). $A x+ B \sin z= C,\qquad A,B,C={\rm const}$;\\
c). $x \cos z+ y \sin z=0$.
}\end{Ex}

\subsection{General structure of sub-Riemannian minimal surfaces and the method of characteristics}
 
Equation (\ref{mineqF}) is essentially degenerate and quite difficult to treat by direct methods.
Instead here we propose an alternative way to study its solutions using the classical method of characteristics.
The first step in this direction is to pass from the sub-Riemannian normal $\nu$ to the
its $\Delta$-orthogonal compliment. Namely, we denote by $\eta$ the horizontal unite vector field $\eta=\eta_1 X_1+\eta_2 X_2\in\Delta$ such that $\eta_1=\nu_2$ and $\eta_2=-\nu_1$. Clearly $\eta$ is well defined on $W\setminus\Sigma$. 
One can easily check that  $\la\nu, \eta\ra_\Delta=0$ and $\eta(q)\in T_qW$ for all $q\in W\setminus\Sigma$. 
The vector field $\eta$ is called the {\it characteristic vector field of $W$}, its integral curves  on $W$ are the {\it characteristic curves of $W$}. By definition, these curves  are horizontal curves $t\mapsto \gamma(t)\in W$
satisfying
$$
\dot\gamma(t)\in T_{\gamma(t)}W\cap\Delta_{\gamma(t)}\qquad \forall t.
$$

Since $\|\nu\|_{\Delta}=\|\eta\|_{\Delta}=1$,  we can introduce a new parameter $\varphi\in\real$ such that
$$
\cos\varphi=\eta_1,\qquad \sin\varphi=\eta_2, 
$$
so that
$$
\eta \equiv \eta e^{\varphi}=\cos \varphi X_1+\sin\varphi X_2.
$$
Here the  upper-index $\varphi$ stresses out the dependence of the vector field $\eta$ on $\varphi$.

Assume $W$ is a small piece of a smooth sub-Riemannian minimal surface without characteristic points.  
Then $\nu$  and  $\eta^\vp$  are well defined, and (\ref{mineq1}) becomes 
\BE\label{mineq3}
-\sin\varphi\, X_2 \varphi-\cos\varphi\, X_1\varphi=\cos\varphi\, c_{12}^1+\sin\varphi\, c_{12}^2.
\EE
The equation above is a quasilinear PDE and we can apply the classical method of characteristics to find its solutions.
Indeed, denote by $q_i$, $i=1,2,3$ some local coordinates on $M$ and let $t\mapsto (q_1(t), q_2(t), q_3(t))$ be 
a smooth (at least $C^1$) curve. Along this curve  $\dot \varphi=\sum\limits_{i=1}^3\frac{\dd \varphi}{\dd q_i}\dot q_i$ with $\dot{\;\;}=\frac{d}{dt}$. Substituting this expression into (\ref{mineq3}) we get the following system of ODE: 
\BE\label{sys}
\left\{\begin{array}{ccl}
\dot q&=&\eta^{\varphi}(q)\\
\dot \varphi&=&-\cos\varphi\, c_{12}^1(q)-\sin\varphi\, c_{12}^2(q)
\end{array}\right..
\EE
By construction, this system is equivalent to (\ref{mineq3}) at non-characteristic points.

\vspace{5 pt}

The described construction motivates of the following geometric interpretation of the problem. 
Denote by $\ov M=\{\bar q=(q,\vp):\,q\in M,\;\vp\in\real\}$. There is an isomorphism $\ov M\simeq S_{\Delta}M$, where $S_{\Delta}M$ is a horizontal spherical bundle over $M$:
$$
S_{\Delta}M=\{(q,v):\;q\in M,\,v\in \Delta_q,\;\|v\|_{\Delta}=1\}.
$$
By $\pi$ we denote the canonical projection  $\ov M\mapsto M$ and set $X_4\equiv\dd_\vp$. Clearly, the Riemannian structure on $\ov M$ associated with the orthonormal frame $\{X_i\}_{i=1}^4$ is compatible with the sub-Riemannian structure on $\Delta$. Since  $[X_i,X_4]\equiv 0$  for $i=1,2,3$, the non-zero structural constants of the extended frame are the same as the ones of the canonical frame $\{X_i\}_{i=1}^3$.

Now consider the {\it generalized characteristic vector field} $V\in Vec(\ov M)$:
$$
V=\cos\vp X_1+\sin \vp X_2 +g X_4, \qquad  g=-c_{12}^1\cos\vp-c_{12}^2\sin\vp.
$$
It is easy to see that $V$ is well defined everywhere on $\ov M$, it has no singular points, and it is not difficult to verify that it is invariant w.r.t. the choice of the horizontal basis on $\Delta$. Moreover, since  
$\pi_*[V(\bar q)]=\eta^{\vp}(q)$ the projection of the integral curves of $V$ on $M$ are exactly the characteristics of equation (\ref{mineq3}).

Let us take now a smooth vector field $\xi\in Vec(\ov M)$ and denote by $\Gamma$ its integral curve starting at 
$\bar q_0$:\;\;$\Gamma(s)=e^{s \xi}\bar q_0$. Further, denote    
\BE\label{gensurf0}
{\mathcal W}=\{\bar q(t,s)=\left(e^{t V}\circ e^{s\xi}\right)\bar q_0, \quad t\in[-\varepsilon_1, \varepsilon_2], \quad s\in[-\delta_1, \delta_2]\},
\EE
where $\delta_i$ and $\varepsilon_i$, $i=1,2$, are some positive numbers, possibly small. By construction, ${\mathcal W}$ 
is a small piece of the solution of the Cauchy problem
\BE\label{Cauchy0}
\left\{\begin{array}{l}
\dot{ \bar q}(t,s)=V(\bar q(t,s))\\
\bar q(0,s)=\Gamma(s),\quad s\in[-\delta_1,\delta_2]
\end{array}\right..
\EE
The theorem of existence and uniqueness of solutions of ODE guarantees that this solution is unique, provided 
the manifold $M$ and the curve of initial conditions  $\Gamma$ are sufficiently smooth. Moreover, if $\xi\wedge V\big|_\Gamma\ne 0$, then $\mathcal W$ locally has a structure of a $2$-dimensional sub-manifold of $\ov M$.

\begin{proposition} The point $\bar q\in\ov M$ is a singular point of ${\mathcal W}$ iff
any vector $\zeta=\sum\limits_{i=1}^4\zeta_i X_i(\bar q)\in T_{\bar q} {\mathcal W}$ satisfy
\BE\label{condA}
\zeta_1 \sin\vp-\zeta_2\cos \vp=0,\quad \zeta_3=0,
\EE
\BE\label{condB}
g(\bar q)\zeta_1=\zeta_4\cos\vp,\qquad g(\bar q)\zeta_2=\zeta_4\sin\vp.
\EE
\end{proposition}
\noindent{\bf Proof}. The point $\bar q\in\ov M$ is a singular point of ${\mathcal W}$ iff
$$
{\rm dim\;span}\{V(\bar q),\zeta\}<2. 
$$
In other words, all second order minors of the matrix 
$$
\left(\begin{array}{cccc}
\zeta_1&\zeta_2&\zeta_3&\zeta_4\\
\cv&\sv&0&g(\bar q)\\
\end{array}\right)
$$
are zero. Taking into account that the pair of conditions $\zeta_3\sv=0$, $\zeta_3\cv=0$ imply $\zeta_3=0$ we 
get (\ref{condA}) and (\ref{condB}). $\square$
  
\vspace{5 pt}

Conditions (\ref{condA}) and  (\ref{condB}) admit a very clear geometrical interpretation. Indeed, needless to say that any singular point of $\mathcal W$ is a singular point of its projection $W=\pi[{\mathcal W}]$. Further, it is easy to see that if non of  conditions (\ref{condA}), (\ref{condB}), but $\zeta_3=0$, is verified, then at $q=\pi[\bar q]$ we have
$\pi_*[T_{\bar q}{\mathcal W}]=\Delta_q$ and hence  $q$ is a characteristic point of $W$. Such a point  is called a {\it regular characteristic point} in he sequel. If at $\bar q$  both conditions (\ref{condA}) are satisfied, then $T_{\bar q}{\mathcal W}={\rm span}\{V(\bar q),\dd_\vp\}$. In this case for all $\zeta\in T_{\bar q}{\mathcal W}$ ${\rm dim\,rank}\{\eta^\vp(q), \pi_*[\zeta]\}=1$, and hence the characteristic point $q$  is a singular point of $W$ (we will say that it is a {\it singular characteristic point}). Finally, if all conditions (\ref{condA}), (\ref{condB}) fail at $\bar q$, then it is a regular point of ${\mathcal W}$, and its projection $q=\pi[\bar q]$ is a regular point of the projected surfaces $W$. Moreover,  by construction,  $W$ satisfies equation (\ref{mineq1}) at $q$. 

Summing up we see that the projected surface $W=\pi[{\mathcal W}]$ is a  $\Delta$-minimal surface, possibly with characteristic points. Moreover, by a suitable choice of the initial curve $\Gamma$ any given piece of a $\Delta$-minimal surface can be presented as the projection of a piece of form (\ref{gensurf0}). Thus in the case of $(2,3)$ contact sub-Riemannian structures Definition \ref{MinS} can be naturally generalized as follows:

\begin{defi}\label{newdef} Let $M$, ${\rm dim} M=3$, be a smooth contact sub-Riemannian manifold. 
We say that the hypersurface $W\in M$ is $\Delta$-minimal w.r.t. the sub-Riemannian structure $\Delta$ of co-rank $1$ iff any piece of it can be presented as the projection of the one-parametric family of solutions of the Cauchy problem (\ref{Cauchy0}) for some curve $\Gamma\in S_{\Delta}M$.
\end{defi}
In what follows we will call ${\mathcal W}$  and $\Gamma$ the {\it generating surface} and {\it generating curve}
of the minimal surface $W=\pi[{\mathcal W}]$.

\vspace{5 pt}

\noindent{\bf Remark}. In a particular but important for the applications case of Lie groups  the  described method 
provides an explicit parameterizations of the minimal surfaces. Indeed, recall that in the Lie group case $c_{12}^1$ and $c_{12}^2$ are constants. Then, for any fixed $s$  we can perform the direct integration of the  second equation of (\ref{sys}). The obtained function $\varphi(t,s)$ can be used to solve the first equation of (\ref{sys}).  If, moreover, 
\BE\label{ruledcond}
c_{12}^1=c_{12}^2=0,
\EE
then $\varphi$ is constant along any characteristic curve. The resulting minimal surface is a kind of ruled surface, whose rulings are the characteristic curves, which are not straight lines in general.

\begin{Ex}{\rm  Consider the case of the Heisenberg group ${\mathbb H}^1$.  Let us use the standard Cartesian coordinates $(x,y,z)$ in $\real^3$ so that the horizontal basis is given by   
$$
X_1=\dd_x-\frac{y}{2}\dd_z,\qquad X_2=\dd_y+\frac{x}{2}\dd_z. 
$$
We fix the orientation by choosing the Reeb vector field $X_3=\dd_z$, so that the only non-zero structural constant is
$c_{12}^3=-1$. Condition (\ref{ruledcond}) is  satisfied, and hence the parameter $\vp$ is constant along characteristic curves. The characteristic vector field reads 
$$
V=\cos\vp\,\dd_x+\sin\vp\,\dd_y+\frac{1}{2}(x\sin \varphi-y\cos\varphi)\dd_z,\qquad\varphi\in[0,2 \pi].
$$
Thus any characteristic curve satisfies the following system of ODE for some fixed $\varphi$: 
\BE\label{H1sys}
\dot x=\cos\vp,\qquad \dot y=\sin\vp,\qquad \dot z=\frac{1}{2}(x\sin \vp-y\cos\vp).
\EE
We easily see that characteristic curves are straight lines.  The minimal surface generated by the curve 
$\Gamma(s)=(x_0(s), y_0(s),z_0(s),\vp(s))$ admits the following parametrization:
\BE\label{HMin}
x(t,s)=t\cos\vp(s)+x_0(s),\qquad y(t,s)=t\sin\vp(s)+y_0(s),
\EE
$$
z(t,s)=\frac{t}{2}(x_0(s)\sin \vp(s)-y_0(s)\cos\vp(s))+z_0(s).
$$
}\end{Ex}

\begin{Ex}{\rm In the case of the group of roto-translations $E_2$ condition (\ref{ruledcond}) is satisfied as well, and 
$\varphi$ is constant along characteristics. The characteristic vector field is given by $V=\cos \varphi \cos z \dd_x+\cos\varphi\sin z \dd_y+\sin\varphi\dd_z$. The characteristic curves are then the integral curves of the system:
\BE\label{E2sys}
\dot x=\cos\vp\cos z\qquad \dot y=\cos\vp\sin z,\qquad \dot z=\sin \vp.
\EE
The solution generated by the curve $\Gamma(s)=(x_0(s),y_0(s),z_0(s),\vp(s))$  has the form 
$$
x(t,s)=\cos (t\sin\vp(s)+z_0(s))\cot\vp(s)+x_0(s),$$
$$
y(t,s)=-\sin(t\sin\vp(s)+z_0(s))\cot\vp(s)+y_0(s),
$$
$$
z(t,s)=t \sin\vp(s)+z_0(s) 
$$
for all $s$ such that $\vp(s)\ne 0,\pi$, and
$$
x(t,s)=\pm t\cos z_0(s)+x_0(s),\quad y(t,s)=\pm t\sin z_0(s)+y_0(s),\quad z(t,s)=z_0(s)
$$
for $s$ where $\varphi(s)= 0$ mod $\pi$. It is not difficult to verify that this surface is smooth w.r.t. $s$ 
if the generating curve is so. Observe that the characteristic curves, except those that correspond to $\vp(s)=0\,{\rm mod}\,\pi$, are not straight lines. 
}\end{Ex}

\subsection{Local structure of singular sets of sub-Riemannian minimal surfaces}

The lifting of the problem to $\ov M$ gives us the key to the complete understanding of the local structure of the sub-Riemannian minimal surfaces, in particular, in the vicinity of their characteristic points. Indeed, let us denote by $\xi^t=\sum\limits_{i=1}^4\xi_i^t X_i=e^{t V}_*\xi$ the push-forward of the vector field $\xi$ by the characteristic flow $e^{t V}$. Consider a small piece of the parametrized surface $\mathcal W\in\ov M$ of form  (\ref{gensurf0}). 
As
$$\frac{d}{ds}\bar q(t,s)=e^{t V}_*\xi(\bar q(t,s)),
$$
it follows that $\xi^t\in T{\mathcal W}$ for all $t$.  
Moreover, for any curve $\beta(s)=\bar q(t(s),s)\in {\mathcal W}$ 
\BE\label{beta}
\frac{\dd\beta(s)}{\dd s}=t'(s)V(\beta(s))+\xi^{t(s)}(\beta(s)).
\EE
If  $\pi[\beta]$ contains a singular point of $W$ at some $\hat s$, then it is tangent to the characteristic vector field at this point.

Our further analysis is based on the following Taylor's expansion of the vector field  $\xi^t$. 
\begin{proposition}\label{PrN}
Let $\xi\in Vec(\ov M)$ and  consider the curve $\bar q_t=e^{t V}\bar q$  starting at some point $\bar q\in\ov M$. Then 
$$
\begin{array}{cr}
\xi_1^t(\bar q_t)=&\xi_1(\bar q)-t(c_{12}^1(\xi_1\sin\vp-\xi_2\cos\vp)-\xi_3(c_{13}^1\cos\vp+c_{23}^1\sin\vp)
+\xi_4\sin\vp)(\bar q)+o(t^2),
\end{array}
$$
\BE\label{Nform}
\begin{array}{cr}
\xi_2^t(\bar q_t)=&\xi_2(\bar q)-t(c_{12}^2(\xi_1\sin\vp-\xi_2\cos\vp)-\xi_3(c_{13}^2\cos\vp+c_{23}^2\sin\vp)-\xi_4\cos\vp)(\bar q)
+o(t^2),\end{array}
\EE
$$
\begin{array}{cr}
\xi^t_3(\bar q_t)=&\xi_3(\bar q)+t(\xi_1\sin\vp-\xi_2\cos\vp)(\bar q)-\frac{t^2}{2}(\xi_4+c_{12}^1\xi_1
+c_{12}^2\xi_2 +o(t))(\bar q_t)\end{array}.
$$
\end{proposition}
\noindent{\bf Proof} The proof of the proposition consists in a straightforward computation. 
We give just a sketch of it for the third component, the remaining formulae can be derived in the same way. 

First of all  we recall that the push-forward operator admits the following representation (see \cite{AgrBook},~\cite{AgrExp}):
\BE\label{ad}
e^{t V}_*\xi=e^{- t \ad V}\xi=\big(Id-t [V,\cdot]+\frac{t^2}{2!}[V,[V,\cdot]]+\dots\big)\xi.
\EE
Let us compute explicitly the first terms of this expansion. We have
$$
[V,\xi]=(V \xi_1+c_{12}^1(\xi_1\sin\vp-\xi_2\cos\vp)-\xi_3(c_{13}^1\cos\vp+c_{23}^1\sin\vp)+\xi_4\sin\vp)X_1+
$$
$$
+(V \xi_2+c_{12}^2(\xi_1\sin\vp-\xi_2\cos\vp)-\xi_3(c_{13}^2\cos\vp+c_{23}^2\sin\vp)-\xi_4\cos\vp)X_2+
$$
$$
+(V\xi_3-(\xi_1\sin\vp-\xi_2\cos\vp))X_3+(V\xi_4-\xi g) X_4.
$$
This expression can be used in order to calculate the second order brackets. 
In particular, after all necessary simplifications for the third component  we obtain
$$
[V,[V,\xi]]_3=V^2\xi_3-2 V(\xi_1\sin\vp-\xi_2\cos\vp)-(\xi_4+c_{12}^1\xi_1+c_{12}^2\xi_2).
$$
Taking into account that for any function $f$ 
$$
e^{-t V}f=f-tVf+\frac{t}{2!}V^2 f+o(t^3),
$$
we get 
$$
\xi^t_3=e^{- t V}\xi_3+t e^{-t V}(\xi_1\sin\vp-\xi_2\cos\vp)-\frac{t^2}{2}(\xi_4+c_{12}^1\xi_1+c_{12}^2\xi_2)+o(t^3).
$$
Now in order to conclude the proof of the formula it is enough to recall that diffeomorphisms acts on functions as the change of variables, i.e., 
$$
e^{-tV}f (\bar q_t)=f(e^{-t V}\bar q_t)=f(\bar q).
$$
The formulae for $\xi_1^t$ and $\xi_2^t$ can be found in the same way. $\square$

\vspace{5 pt}

Let $\bar q$ be a regular point of $\mathcal W$. Without loss of generality we can assume $\bar q=\bar q(0,0)$. 
Since $\bar q$ is  regular, there exists a vector $\zeta\in T_{\bar q}{\mathcal W}$  such that $\zeta\wedge V(\bar q)\ne 0$.  Let $\xi\in T{\mathcal W}$ be a smooth vector field such that $\xi(\bar q)=\zeta$. 
We can take as a generating curves a small piece of the curve $\bar q(s)=e^{s \xi}\bar q$, $s\in[-\delta_1, \delta_2]$, passing through $\bar q$. The characteristic points of $W=\pi[{\mathcal W}]$ correspond to zeros of the function
$$
\phi(t,s)=\xi^t_3(e^{t V}\bar q(s))=\phi_0(s)+t\phi_1(s)-\frac{t^2}{2}\phi_2(t,s), 
$$
where $\phi_0(s)=\xi_3(\bar q(s))$, $\phi_1(s)=(\xi_1\sin\vp-\xi_2\cos\vp )(\bar q(s))$, and
$\phi_2(t,s)=(\xi_4+c_{12}^1\xi_1+c_{12}^2\xi_2+t f)(\bar q(t,s))$. Here the function $f$ contains the higher order terms of the expansion (\ref{Nform}).

The point $q=\pi[\bar q]$ is  non-characteristic if and only if  $\phi_0(0)\ne 0$ for some choice of $\zeta$. Now let us see what happens if $q$ is a characteristic point of $W$.

\vspace{5pt}

\noindent{\bf Case 1: regular characteristic point}. In this case $\phi_0(0)=0$ for all $\zeta\in T_{\bar q}{\mathcal W}$ and it is possible to choose  $\zeta$ in such a way that  $\phi_1(0)\ne 0$. Therefore we can write $\phi(t,s)=\phi_0(s)+t\hat\phi_1(t,s)$, where $\hat \phi_1=\phi_1-\frac{t}{2}\phi_2$. We have
$$
\phi(0,0)=0, \qquad \frac{\dd\phi}{\dd t}(0,0)=\phi_1(0)\ne 0.
$$
By the implicit function theorem in a small neighborhood of the origin in $\real ^2$ there exists a unique curve $t=t(s)$ such that $t(0)=0$ and $\phi(\bar q(t(s),s))=0$. By construction, the curve $q(t(s),s)\in W$ consists of regular characteristic points, we will call such a curve a {\it simple singular curve} of $W$.\footnote{In the language of the singularity theory the simple singular curves are resolvable singularities of sub-Riemannian minimal surfaces.}

\vspace{5 pt}

\noindent{\bf Case 2: generic singular  characteristic point}. In this case 
$\phi_0(0)=0$ and $\phi_1(0)=0$ for all $\zeta\in T_{\bar q}{\mathcal W}$.
Thus we can assume that $\phi_0(s)=s\tilde\phi_0(s)$ and $\phi_1(s)=s\tilde\phi_1(s)$ for some functions $\tilde \phi_0$ and $\tilde \phi_1$. Observe that $\tilde\phi_0(0)=\frac{\dd}{\dd s}\xi_3(\bar q(s))\big|_{s=0}$.
First let us consider in detail the generic situation when $\tilde\phi_0(0)\ne 0$. 
\begin{Theorem}\label{Classification} Let $\bar q$ be a regular point of a  generating surface ${\mathcal W}$ such that $q=\pi[\bar q]$ is a  generic singular  characteristic point of $W=\pi[{\mathcal W}]$. Then:\\
a) a small enough  neighborhood of $q$ in $W$ contains a pair of simple singular curves, which touch each other at $q$ and form a unique smooth curve on $W$ passing through $q$;\\
b) there exists a choice of local coordinates $(\tau,\sigma)\in\Omega\subseteq \real^2$ such that   
$$
W=\{q(\tau,\sigma)=(\tau^2,\sigma,\tau\sigma),\;\; \tau, \sigma\in \Omega\subseteq\real^2\},\qquad q=q(0,0), 
$$
i.e., in the neighborhood of $q$ $W$ has the structure of Whitney's umbrella.
\end{Theorem}
\noindent{\bf Proof}. First of all we observe that since $\bar q$ is a regular point $\phi_2(0,0)\ne 0$. Indeed, otherwise we would have $\zeta_4=-c_{12}^1\zeta_1-c_{12}^2\zeta_2$ and $\zeta_1\sin\vp=\zeta_2\cos\vp$. It is not difficult to see that these two conditions are equivalent to (\ref{condB}), and this contradicts the regularity assumption on the initial point $\bar q$.

In order to prove  part a) let us treat the implicit equation $\phi(t,s)=0$ as a quadratic equation w.r.t. $t$-variable.
The discriminant of this equation is given by the function
$$
D(t,s)=s(s\tilde \phi_1^2(s)+2\tilde\phi_0(s)\phi_2(t,s))=s\beta(t,s), 
$$
where  $D(0,0)=0$ and $\beta(0,0)\ne 0$. We have $\frac{\dd}{\dd t}D(0,0)=0$ and 
$\frac{\dd}{\dd s}D(0,0)=\beta(0,0)\ne 0$. Thus  we can assume $D(0,s)>0$  in a small neighborhood of $(0,0)$ for $s>0$.
Indeed, this condition can be always satisfied by an appropriate choice of the sign of the parameter $s$. Now we obtain two implicit equations
\BE\label{tpm}
t=\frac{s\tilde\phi_1(s)\pm\sqrt{s\beta(t,s)}}{\phi_2(t,s)},
\EE
which describe zero level sets of two functions
$$
\Phi^{\pm}(t,s)=t-\frac{s\tilde\phi_1(s)\pm\sqrt{s\beta(t,s)}}{\phi_2(t,s)}.
$$
We have 
$$
\Phi^{\pm}(0,0)=0,\qquad \frac{\dd}{\dd t}\Phi^{\pm}(0,0)=1.
$$
Applying the implicit function theorem to each of the functions $\Phi^{\pm}$ we see that in a small enough neighborhood of  the origin in the half-plane $s\ge 0$ of the $(t,s)$-plane there exist two curves $(t^{\pm}(s),s)$ satisfying (\ref{tpm}). Since $t^\pm(0)=0$ these curves meet each other at the origin of the $(t,s)$-plane. Since $\lim\limits_{s\to 0}\frac{\dd}{\dd s}\Phi^{\pm}(0,s)=\mp\infty$ they are both tangent to the $t$-axis. The corresponding curves $q(t^{\pm}(s),s)$ on $W$ form a unique curve, which consist of two branches, which smoothly glue together at the singular characteristic point $q(0,0)$. Moreover, taking into account (\ref{beta}) we see that the curve formed by the pair of  curves which bifurcate at a singular characteristic point is tangent to the characteristic passing through this point.

We claim that the small enough pieces of  curves $q(t^{\pm}(s),s)$, $s>0$, contain only regular characteristic points. 
In order prove this it is enough to show that the function $\chi=\xi_1^t\sin\vp-\xi_2^t\cos\vp$ is different from zero in a small neighborhood of $\bar q$. First of all we observe that for $s$ fixed
$$
\sin\vp(t,s)=\sin(\vp(s)+t \xi_4(s)+o(t^2;s))=\sin\vp(s)+t\xi_4(s)\cos\vp(s)+o(t^2;s),
$$
$$
\cos\vp(t,s)=\cos(\vp(s)+t \xi_4(s)+o(t^2;s))=\cos\vp(s)-t\xi_4(s)\sin\vp(s)+o(t^2;s).
$$
Hereafter, for briefness, $\vp(s)=\vp(0,s)$, $\xi_i(s)=\xi_i(\bar q(0,s))$.

Without loss of generality we can assume that $\xi_1(\bar q)=\xi_2(\bar q)=0$.
We have   $\chi(0,0)=0$, and by definition, $\chi(0,s)=\phi_1(s)$. 
Using expressions (\ref{Nform}) we obtain 
$$
\frac{\dd\chi}{\dd t}\Big|_{(0,0)}=\frac{\dd}{\dd t}(\xi^t_1\sin\vp-\xi_2^t\cos\vp)(\bar q(0,0))=
$$
$$
=\zeta_4(\zeta_1\cos\vp(0)+\zeta_2(0)\sin\vp(0)-1)=-\zeta_4\ne 0,
$$ 
because, by assumption, $\bar q$ is a regular point of ${\mathcal W}$.
So if $\phi_1(s)\ne 0$ for $s>0$, then  it follows that function $\chi$ is different from zero in a small neighborhood of $\bar q$, and hence in this neighborhood the points $q(t^{\pm}(s),s)$, $s>0$, are regular characteristic points.

Let us prove part b). The sub-Riemannian minimal surfaces are the images of the map 
$q:\;I\subseteq\real^2\mapsto M$ such that $q(t,s)=\pi[\bar q(t,s)]\in M$, $(t,s)\in I\subseteq\real^2$. 
We may assume again that $W$ is generated by an integral curve of the field $\xi\in \ov M$ with  $\xi(\bar q)=\zeta$ and  $\zeta_i=0$, $i=1,2,3$. In addition, by assumption,  $\frac{\dd}{\dd s}\xi_3(\bar q(s))\big|_{s=0}\ne 0$. We have
\BE\label{first}
\frac{\dd q}{\dd t}\Big|_{(0,0)}=\big(\cos\vp X_1+\sin\vp X_2\big)(q(0,0)),\quad
\frac{\dd q}{\dd s}\Big|_{(0,0)}=\pi_*[\zeta]=0.
\EE
In order to proof b), according to the well-known result by Whitney (\cite{Whit}), it is enough to show that the vectors   $v_1=\frac{\dd q}{\dd t}\Big|_{(0,0)}$, $v_2=\frac{\dd^2 q}{\dd t\dd s}\Big|_{(0,0)}$ and $v_3=\frac{\dd^2 q}{\dd s^2}\Big|_{(0,0)}$ are linearly independent. We have
$$
\frac{\dd q (t,s)}{\dd t\dd s}=\vp'_s(t,s)\big(-\sin\vp(s)X_1(q(t,s))+\cos\vp(t,s)X_2(q(t,s))\big)+
$$
$$
+\cos\vp(t,s)\frac{\dd}{\dd s}(X_1(q(t,s))+\sin\vp(t,s)\frac{\dd}{\dd s}(X_2(q(t,s)).
$$
By (\ref{first}), the last two terms on the expression above vanish at $q(0,0)$. So, 
$$
v_2=\xi_4(0)\big(-\sin\vp X_1+\cos\vp X_2\big)(q(0,0)),
$$
since $\vp'_s(s)=\xi_4(s)$. Analogously we get
$$
v_3=\left(\frac{\dd \xi_1}{\dd s}X_1+
\frac{\dd \xi_2}{\dd s}X_2+\frac{\dd \xi_3}{\dd s}X_3\right)(q(0,0)).
$$
An easy calculation shows that $\det\{v_1,v_2,v_3\}=\xi_4(0)\frac{\dd}{\dd s}\xi_3(0)\ne 0$. $\square$

\vspace{5 pt}

We see that the assumption $\frac{\dd}{\dd s}\xi_3(\bar q(s))\big|_{s=0}\ne 0$ 
implies that a generic isolated singular point of a sub-Riemannian minimal surface  gives rise to a self-intersection
of Whitney's umbrella type. The analysis of non-generic singularities may be done in each single case, but this lies behind  the scope of this paper. We just want to mention the other possible configurations of simple singular curves which appear when  the function $\phi$ has a higher order tangency to zero at the origin. For instance, if $W$ contains a curve, which passes  through the singular point $q(0,0)$ and  such that $\phi_1(s)\equiv 0$, then the same argument as before shows the existence of a pair of simple singular curves touching each other at $q(0,0)$ and having a common tangent (actually they both are tangent to the characteristic passing through $q(0,0)$). In addition, assume that the first non-zero derivative of $\phi_0$ at the origin is of order $k$, $k\ge 2$. Then if $k$ is even, each of the singular curves is regular at $q(0,0)$, while if $k$ is odd each singular curve have a cusp point at $q(0,0)$.  Possible configurations of simple singular curves in the neighborhood of an isolated singular point are shown in Fig.1.
\ppotR{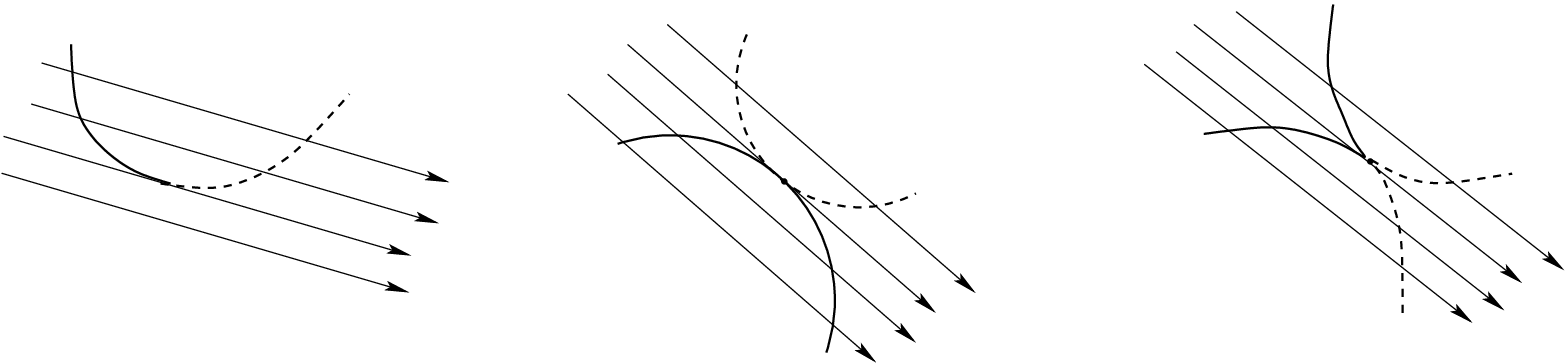}{Simple singular curves in the neighborhood of a singular characteristic point}{2.7}

If $\phi_0(s)\equiv\phi_1(s)\equiv 0$, then ${\mathcal W}$ contains a curve everywhere tangent to the plane 
${\rm span}\{\dd_\vp, V\}$, let us denote it by $\Gamma_*$. Its projection $\gamma_*=\pi[\Gamma_*]$ consists of singular characteristic points ({\it strongly singular curve}). Any narrow enough stripe along 
$\gamma_*$ contains no other characteristic points of $W$. To see this  we can take $\Gamma_*$ as the generating curve.
We can always do that provided the assumption $\xi\wedge V\big|_{\Gamma_*}\ne 0$ is satisfied.
Then at any $s$ the equation $\phi(t,s)=0$   has only trivial solution $t=0$.
Moreover, $\phi(s,t)\ne 0$ for small enough $t>0$ since $\phi_2(t,s)\ne 0$.

\vspace{5pt}

\noindent{\bf Case 3: isolated singular point}. A very degenerate situation occurs when the surface $\mathcal W$ contains a purely "vertical curve" $\bar q(s)=(q,\vp(s))$ passing through $\bar q$. Then the whole strongly singular curve $\gamma_*$  collapses into a single point $q=\pi[\bar q]$.  The same argument as before, after an obvious modification, implies that $q$ is an isolated singular point. Since by assumption $\bar q(0,0)$ is a regular point of ${\mathcal W}$, the forth component $\xi_4(s)=\vp'(s)\ne 0$. Therefore the characteristic vector $\eta^\vp(q(0,0)$ rotates monotonically in the plane 
$T_{q(0,0)}W$. In particular, it follows that the index of the isolated singular point is $+1$. 
  
  \vspace{5pt}

We now conclude this paragraph by the following classification of characteristic points of sub-Riemannian minimal
surfaces:\\
Assume $\bar q$ is a regular point of the surface $\mathcal W\in \ov M$ of form (\ref{gensurf0}) such that 
$q=\pi[\bar q]$ is a characteristic point of the projected surface $W=\pi[{\mathcal W}]$. Then in a small neighborhood 
of $q\in W$ there realizes one of the following situations:
\begin{itemize}
\item $q$ is a regular point of $W$. In a small enough neighborhood of $q$ the surface $W$ contains a unique simple
singular curve passing through $q$;
\item $q$ is a singular point of $W$. In addition,
\begin{itemize}
\item it can be an isolated singular characteristic point;
\item in a small enough neighborhood of $q$ the surface $W$ can contain a unique strongly singular 
curve passing through $q$;
\item in a small enough neighborhood of $q$ the surface $W$ can contain a pair of simple singular curves touching each other at  $q$.
\end{itemize}
\end{itemize}

\subsection{Example: singular sets of sub-Riemannian minimal surfaces in ${\mathbb H}^1$}

Let us illustrate the results of previous subsection  by examples of minimal surfaces associated to the  Heisenberg distribution in ${\mathbb H}^1$. Some of the facts that will be discussed below were already noticed in \cite{Malc} and \cite{Hw}.
We limit ourselves to consider only the sub-Riemannian minimal surfaces generated by smooth
generating curves.  Due to the explicit parameterization (\ref{HMin}) knowing the generating curve $\Gamma(s)$, $s\in\real$, is enough to reconstruct the whole projected surface together with its singular set. 

We start by an observation that since the Heisenberg group is nilpotent, the Lie brackets of the generalized  
characteristic vector field $V$ with the vector fields of the canonical frame  $\{X_i\}_{i=1}^4$ of order greater than $2$ vanish. This implies that for any vector field  $\xi\in Vec(\ov M)$ expansions (\ref{Nform}) contains at most quadratic terms w.r.t. $t$. More precisely, 
$$
\xi^t(q_t)\equiv e^{t V}_* \xi(e^{t V}\bar q)=\sum\limits_{i=1}^4 \alpha_i(t, \bar q) X_i(\bar q_t)
$$
where
$$
\alpha_1(t,\bar q)=\xi_1(\bar q)-t\,\xi_4(\bar q)\sv(\bar q),\qquad \alpha_2(t, \bar q)=\xi_2(\bar q)+t\,\xi_4(\bar q)\cv(\bar q),
$$
\BE\label{alpha}
\alpha_3(t, \bar q)=\xi_3(\bar q)+t(\xi_1(\bar q)\sv(\bar q)-\xi_2(\bar q)\cv(\bar q))-\frac{t^2}{2}\xi_4(\bar q),
\EE
$$
\alpha_4(t, \bar q)=\xi_4(\bar q).
$$
Moreover, $\zeta\wedge V(\bar q)=0$ for $\zeta\in T_{\bar q}M$ if and only if
$$
\zeta_1\sin\vp-\zeta_2\cos\vp=0,\qquad\zeta_3=\zeta_4=0.
$$

\begin{proposition}\label{Pr3} Let  $s\mapsto \Gamma(s)$ be  an integral curve of a smooth vector field $\xi\in Vec(\ov M)$ such that $\xi\ne 0$ and $\xi\wedge V\big|_\Gamma\ne 0$. Let  $W$ be a sub-Riemannian minimal surface generated by $\Gamma$. Then if for some  $\hat s$ 
\BE\label{D0}
2\xi_3\xi_4+(\xi_1 \sin\vp-\xi_2\cos\vp)^2\Big|_{\Gamma(\hat s)}=0,
\EE then the point $q(\hat t,\hat s)$ is a singular point of $W$ for
\BE\label{hatT}
\hat t=\frac{\xi_1 \sin\vp-\xi_2\cos\vp}{\xi_4}\Big|_{\Gamma(\hat s)}.
\EE
Moreover, the singular points of the type $(\hat t, \hat s)$, where $\hat t$ is defined by (\ref{hatT}) are the only  
singular points of $W$. 
\end{proposition} \label{spoint} 
\noindent{\bf Proof.} First of all we observe that, due to the non-degeneracy condition $\xi\wedge V\big|_\Gamma\ne 0$,
if at some $\hat s$  we have (\ref{D0}), then $\xi_4(\hat s)\ne 0$. Further, the point $q_*=\pi[\bar q(\hat t,\hat s)]$ is a singular point of $W$ iff 
$$
{\rm dim\; span}\{\pi_*[\xi^{\hat t}],\eta^\vp\}\Big|_{q_*}<2,
$$
i.e.,
\BE\label{minors3}
\cv(\hat s)\alpha_2(\hat t,\hat s)-\sv(s)\alpha_1(\hat t,\hat s)=(\cv(\hat s)\xi_2(\hat s)-\sv(\hat s)\xi_1(\hat s))+t\xi_4(\hat s)=0,   
\EE
$$
\cv(\hat s)\alpha_3(\hat t,\hat s)=0,\qquad \sv(\hat s)\alpha_3(\hat t,\hat s)=0.
$$
Taking into account (\ref{alpha}), one can easily see that (\ref{D0}) and (\ref{hatT}) imply (\ref{minors3}) and vice versa.
$\square$
\vspace{5 pt}

The left-hand side of $(\ref{D0})$ is the discriminant 
\BE\label{D}
D(s)=(\xi_1(s)\sv(s)-\xi_2(s)\cv(s))^2+2\xi_3(s)\xi_4(s)
\EE
of the quadratic equation
\BE\label{Scurve}
\alpha_3(t,s)=\xi_3(s)+t(\xi_1(s)\sv(s)-\xi_2(s)\cv(s))-\frac{t^2}{2}\xi_4(s)=0,
\EE
which describes the characteristic points of $W$.  If $\xi_4(s)\ne 0$,  then it has at most $2$ real roots:
$$
t_*^{\pm}(s)=\frac{\xi_1(s)\sv(s)-\xi_2(s)\cv(s)\pm\sqrt D}{\xi_4(s)}
$$ 
provided  $D(s)\ge 0$. Thus the two branches of simple singular curves described in Theorem \ref{Classification} have the form 
\BE\label{gI}
x^\pm_*(s)=t_*^\pm(s) \cv(s)+x_0(s),\qquad y^\pm_*(s)=t_*^\pm(s)\sv(s)+y_0(s),
\EE
$$
z^{\pm}_*(s)=\frac{t_*^\pm(s)}{2}(x_0(s)\sv(s)-y_0(s)\cv(s))+z_0(s).
$$
This curves may touch each other at a singular point, which according to Proposition \ref{Pr3}, 
is a root of the equation $D(s)=0$.  

\begin{Ex}\label{ex1}{\rm The curve $\Gamma(s)=(0,0,s,s)$ for $s\in \real$
generates the counter-clockwise helicoid $q(t,s)=(t \cos s,t \sin s, s)$, whose rulings are parallel to the $(x,y)$-plane.
Since  $\xi_1(s)=\xi_2(s)=0$ and  $\xi_3(s)=\xi_4(s)=1$, we found that $D(s)=2$. Hence the singular set consists of two simple singular curves $(\pm \sqrt 2\cos s,\pm \sqrt 2 \sin s, s)$ which never meet each other, see Fig.2. It is easy to see that all counter-clockwise helicoids have the same structure of singular set. Notice, that the clockwise oriented helicoids, for instance,  the one generated by $\tilde\Gamma(s)=(0,0,-s,s)$, contain no singular curves neither singular points. 

All these facts obviously hold true for all helicoids with rulings parallel to any contact  plane 
$\Delta_p$ and generated by the curve $\Gamma(s)=(p+s w,\vp_0\pm s)$, where  $p=(x_0,y_0,z_0)$ and 
$w=\left(\frac{y_0}{2},-\frac{x_0}{2},1\right)$. Indeed, these helicoids can be obtained just  by shifting the origin to the point $p$ in the previous example. In particular, it follows that the helicoids of the described class have no singular points. The similar result was also obtained in \cite{Hw}.
}
\end{Ex}

\begin{Ex}\label{ex2}{\rm Consider the curve   $\Gamma(s)=(0,0,\frac{1}{6}s^3,s)$ for $s\in [0,2\pi]$.
This curve is the integral curve of the vector field
$$
\xi=\frac{3 z^{2/3}}{4}X_3 +X_4
$$
starting at the point $\bar q_0=(0,0,0,0)$. In particular, for any $s\in[0,2\pi]$  we have 
$\xi_1(s)=\xi_2(s)=0$, $\xi_3(s)=\frac{1}{2}s^2$ and $\xi_4(s)=1$, so that $D(s)=s^2$.
Therefore $q_0$ is a unique singular point of the resulting minimal surface.
There are two simple singular curves $g^\pm(s)=(\pm s \cos s,\pm s \sin s,\frac{1}{6}s^3)$, which bifurcate from the point $q_0$, moreover, they touch each other at this point. In Fig. 3 we show the general look of this surface (Fig. 3a) and the structure of its singular set (Fig. 3b).  We also want to notice that the generating curve of this example does not satisfy the genericity assumption  and the singular point $q_0$ does not give rise to a self-intersection. We will give a simple criterion for the existence of self-intersections for the minimal surfaces in ${\mathbb H}^1$ at the end of this section.   
}
\end{Ex}

\begin{Ex}\label{ex3}{\rm The surface generated by the curve
$$
\Gamma(s)=(\sin s,-\cos s,1+\frac{s}{2},s).
$$
contains a strongly singular curve, actually is it the curve $\gamma=\pi[\Gamma]$. 
Indeed, $\Gamma$ is the integral curve of the field  $\xi=\cos \vp X_1+\sin \vp X_2 +\dd_{\vp}$. 
Notice, that despite $\pi_*[\xi]=\eta^\vp$ we still have $V\wedge\xi\big|_{\Gamma}\ne 0$ provided the forth component
of the field $\xi$ is different from zero. One can easily check that $D\equiv 0$, $t_*(s)\equiv 0$,  and the 
the strongly singular curve is tangent to the characteristic field at every point.
}\end{Ex}

\ppotRR{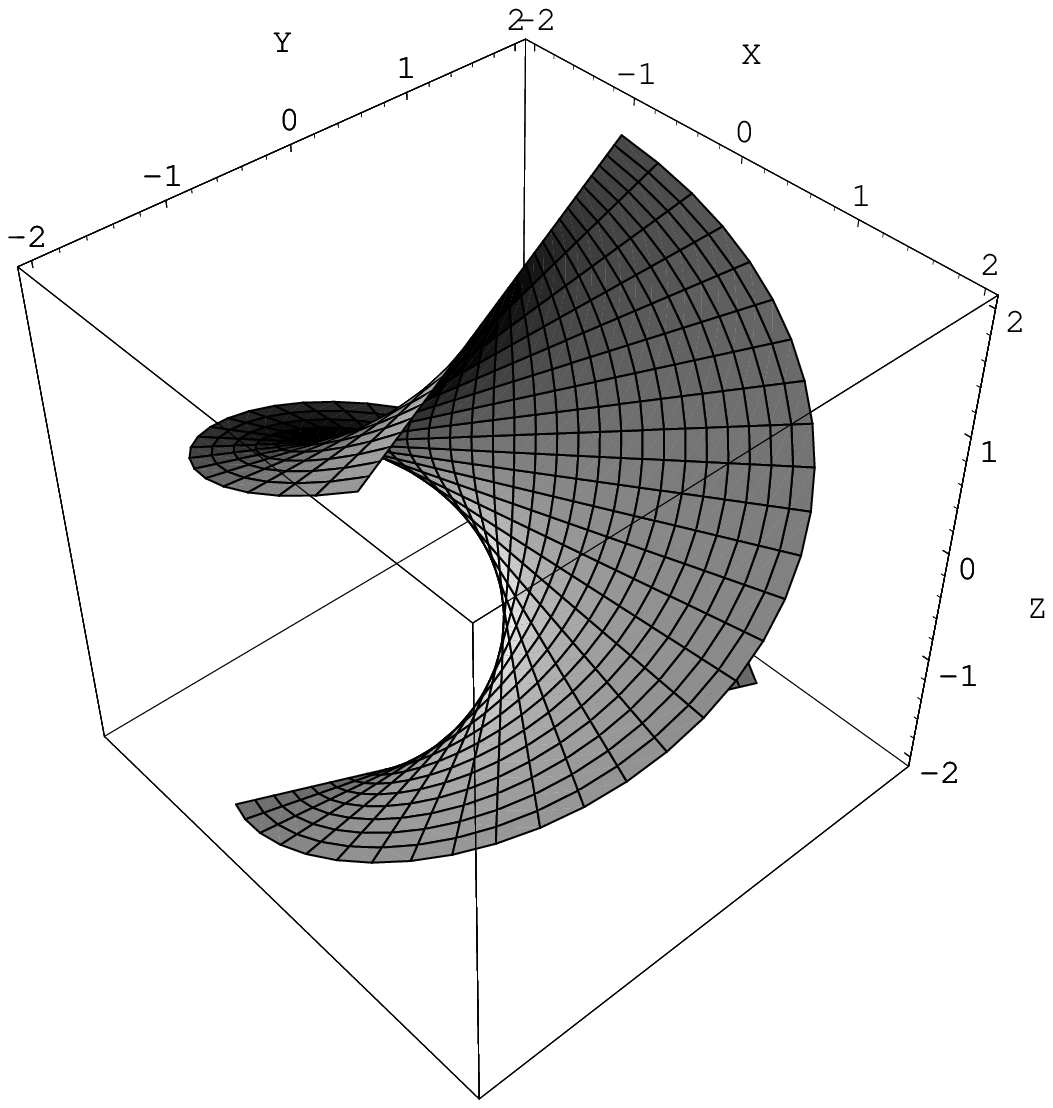}{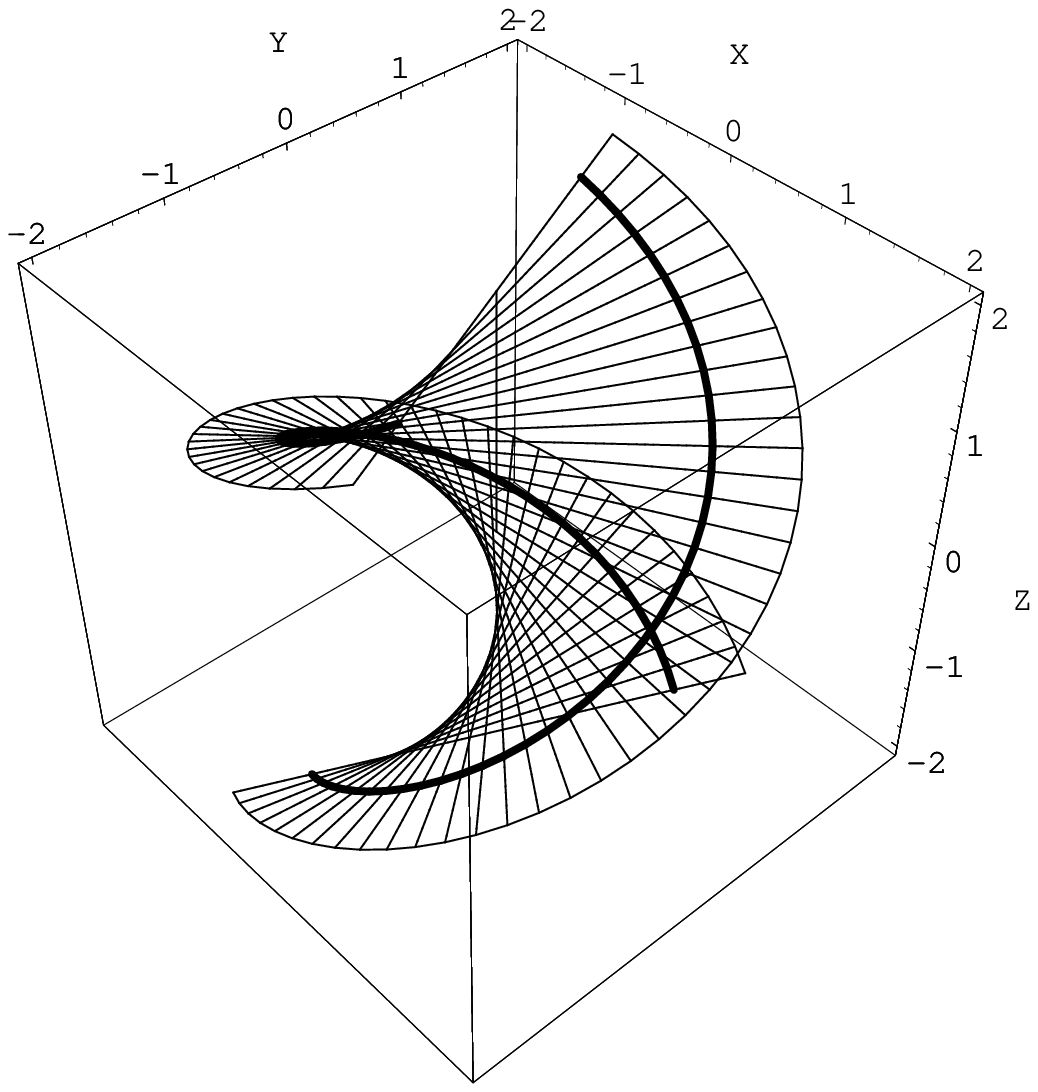}{Counter-clockwise oriented helicoid $(t \cos s,t \sin s,s)$ and its singular set}{6.0}{6.0}
\ppotRR{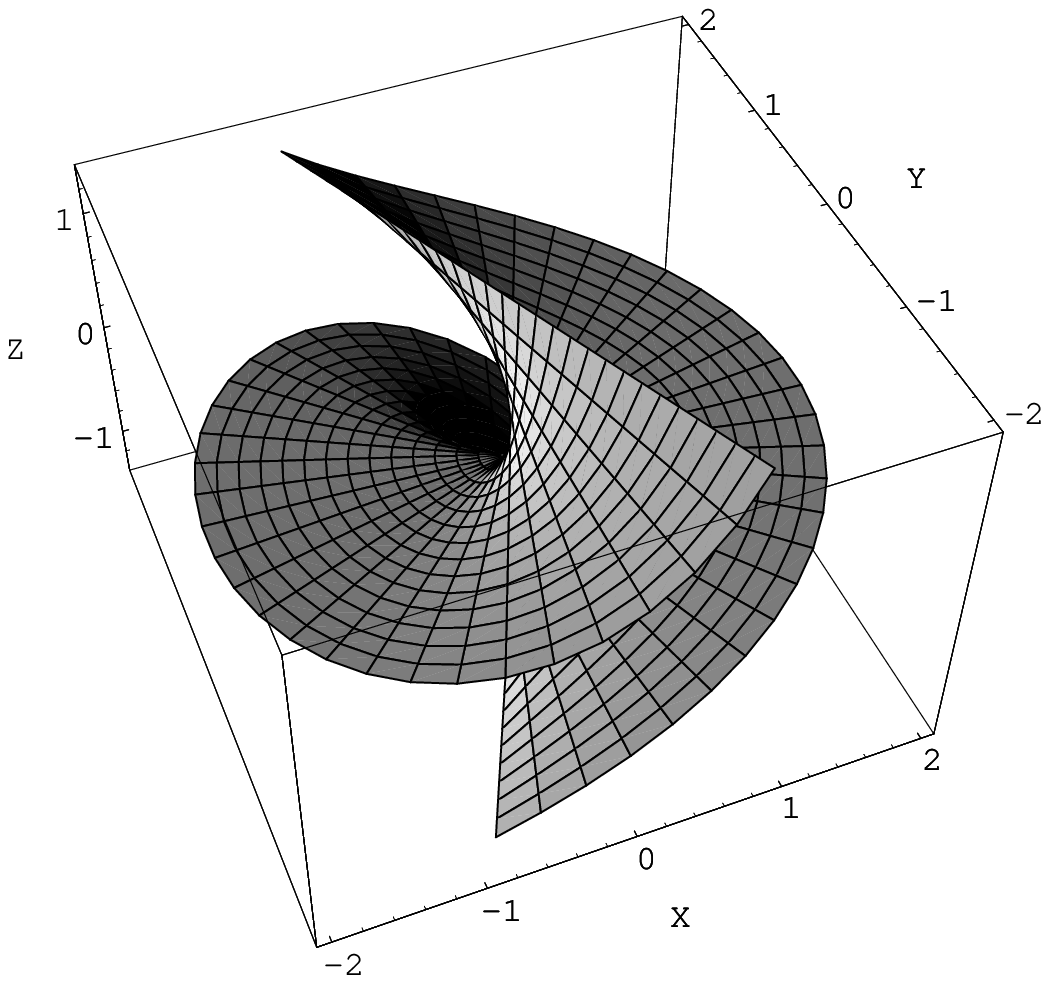}{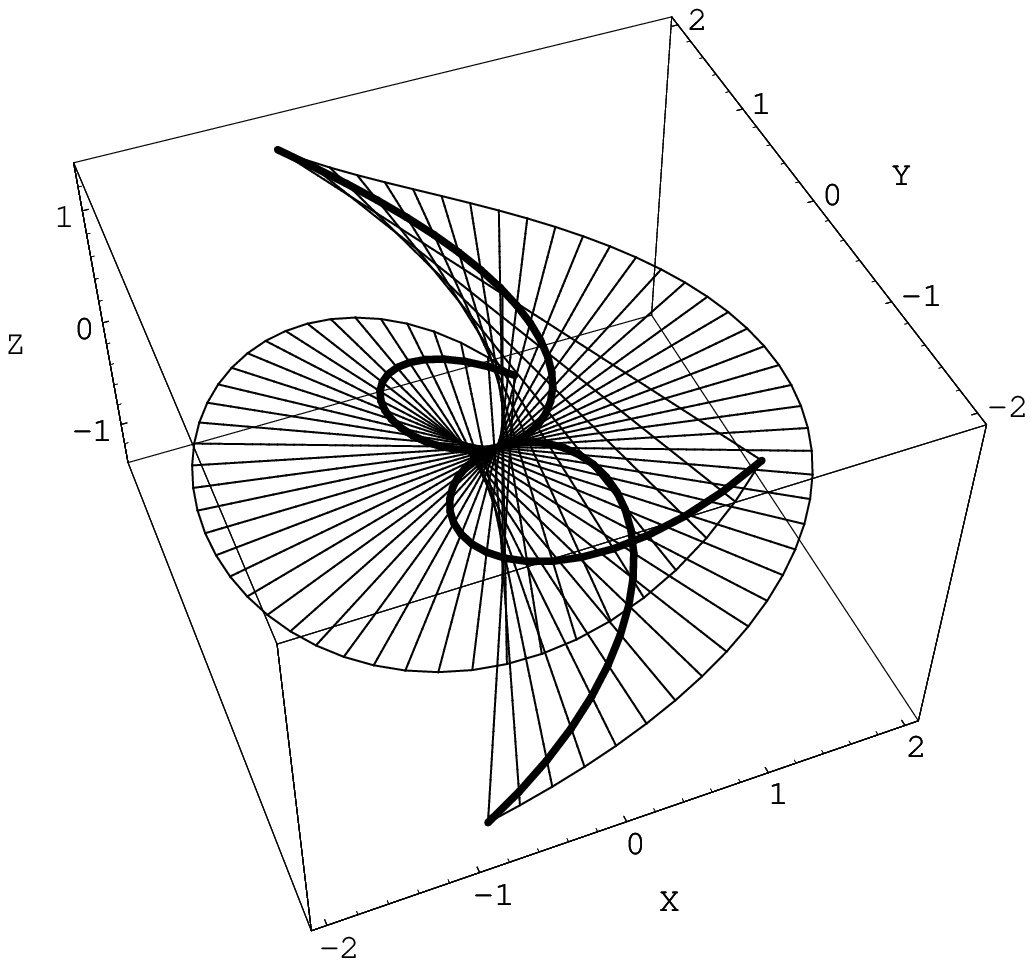}{An example of a surface with  a pair of simple singular curves touching each other at a unique singular point}{6.0}{6.0}
\ppotRR{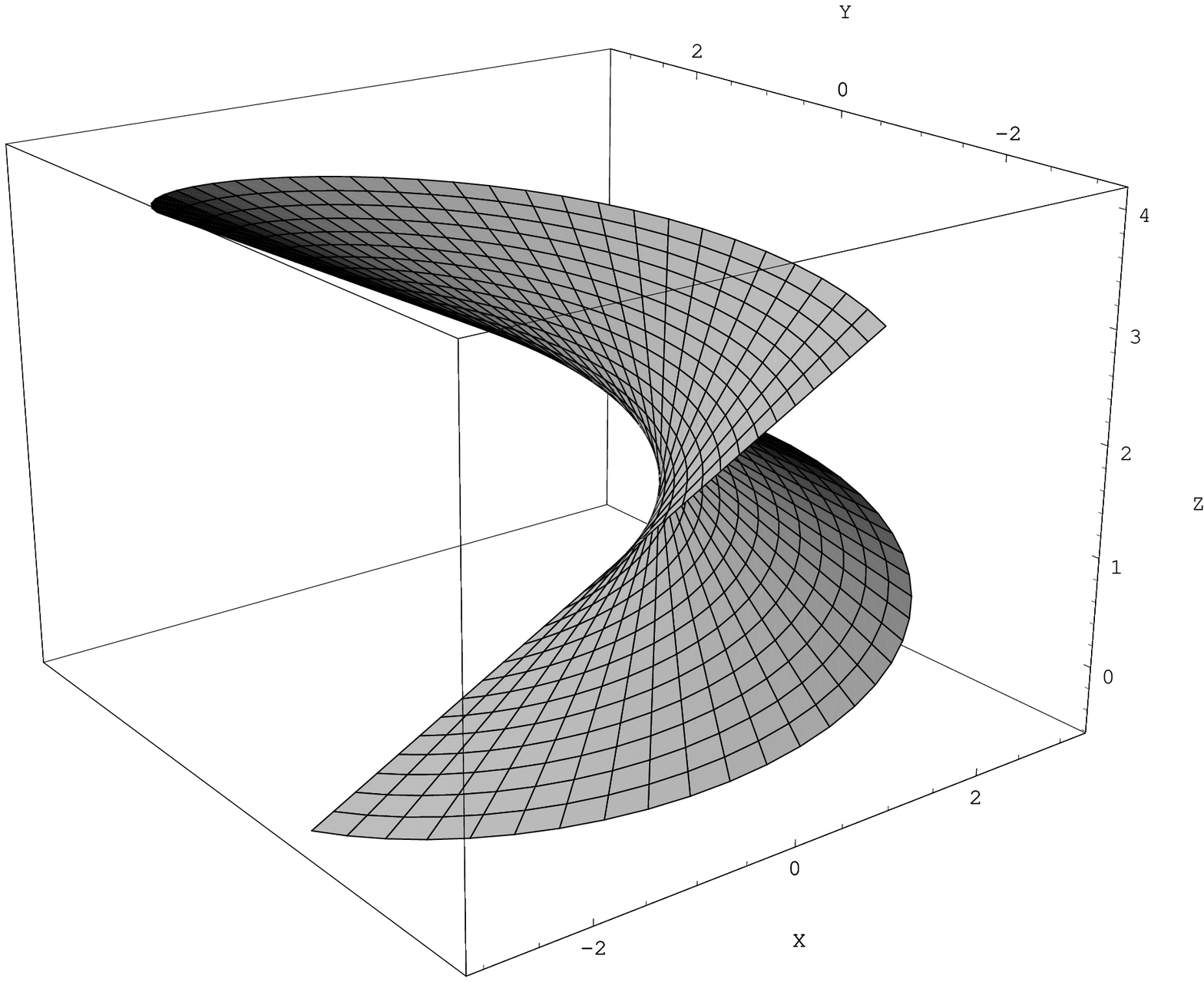}{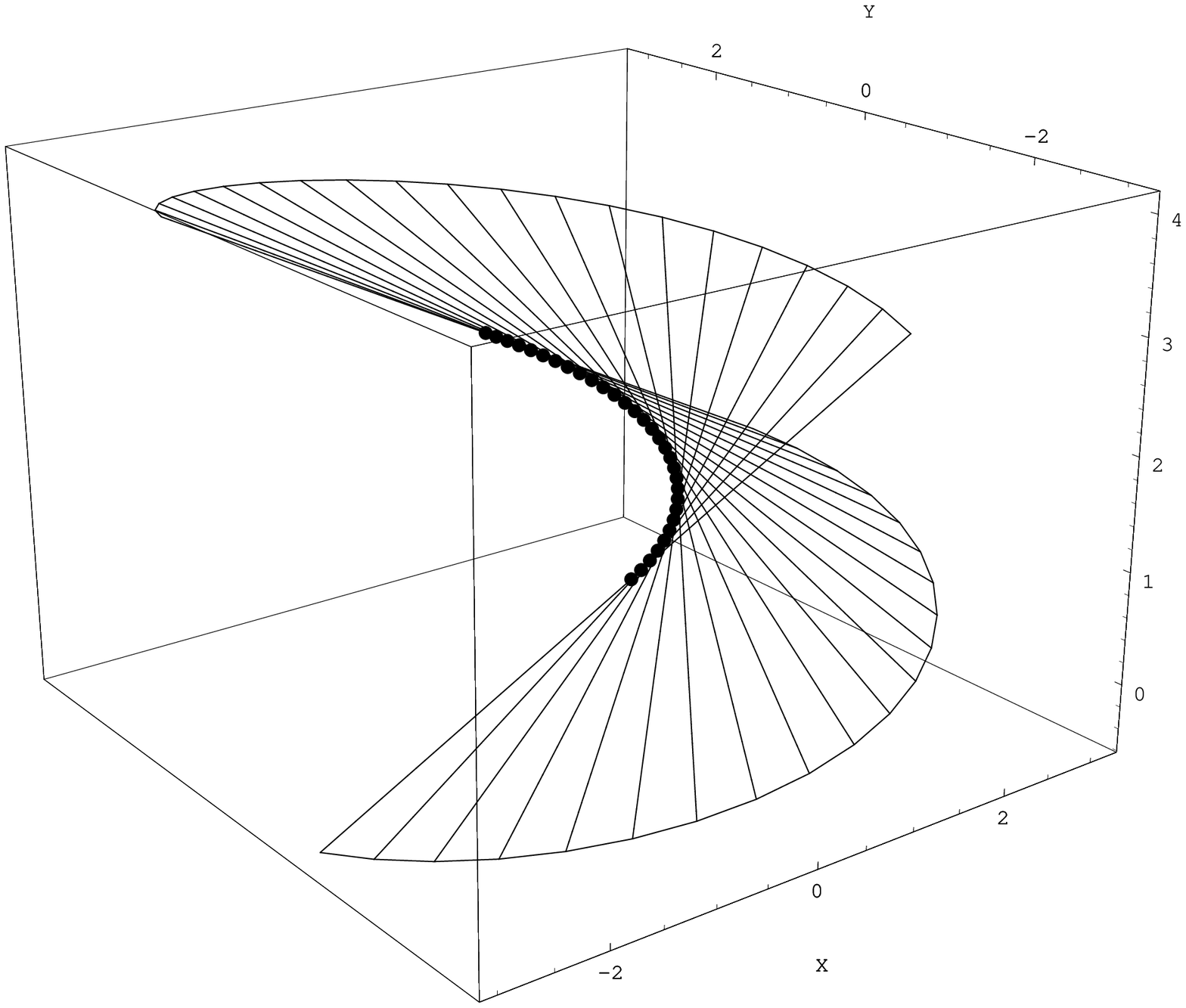}{A piece of a surface containing a strongly singular curve}{6.0}{6.0}

As we know, projections of generating surfaces which contain a purely vertical line $\Gamma_\vp(s)=(x_0,y_0,z_0,s)$ have isolated singular points. Since in the Heisenberg case the characteristic curves are straight lines, the smooth sub-Riemannian minimal surface can contain at most one isolated singular point (this facts was first noticed in \cite{Malc}). Taking $\Gamma_\vp$ as the generating curve we see that the resulting minimal surface is formed by a one-parametric family of ellipses $(t\cos\vp+x_0, t\sin\vp+y_0,\frac{t}{2}(x_0 \sin\vp-y_0 \cos\vp)+z_0)$, $t\ge 0$, that fills the whole plane 
$\Delta_p$, $p=(x_0,y_0,z_0)$. In particular, it follows that in ${\mathbb H}^1$ the only the sub-Riemannian minimal surfaces having isolated characteristic points are planes.
\vspace{5 pt}

We have already seen that any generic singular characteristic point is a starting point of a germ of a curve of self-intersections of Whitney's umbrella type. The next proposition describes all possible self-intersections for minimal surfaces in ${\mathbb H}^1$ case.  This result is an immediate consequence of the explicit parametrization (\ref{HMin}). Here we denote  $v(s)=\eta^\vp(s)(\gamma(s))\in \real^3$.

\begin{proposition}\label{Pr4}
Let $W=\{q(t,s):\; s\in[0,s_1],\, t\in \real\}$ be a piece of sub-Riemannian minimal surface in ${\mathbb H}^1$
corresponding to the generating curve $\Gamma(s)=(\gamma(s), \vp(s))$.
If there exist a pair $a,\,b\in [\delta_1,\delta_2]$ such that $a\ne b$  and a pair of numbers $\tau_1,\,\tau_2\in \real$ such that
\BE\label{self}
\gamma(a)-\gamma(b)=\tau_1 v(a)+ \tau_2 v(b),
\EE
then $W$ contains a point of self-intersection $q_*=q(a,-\tau_1)=q(b, \tau_2)$.
\end{proposition}

\begin{Ex}\label{ex4}{\rm Let $\Gamma(s)=(-2 \cos s,-2 \sin s,\cos s,s)$ be the generating curve and denote
$\gamma=\pi[\Gamma]$. We have $\xi_1(s)=2 \sin s$, $\xi_2(s)=-2 \cos s$, $\xi_3(s)= -2-\sin s$ and $\xi_4(s)=1$. First
we find the singular points. In the present case  $D(s)=-2 \sin s$, i.e., according to Proposition \ref{Pr3}, the sub-Riemannian minimal surface generated by $\Gamma$ contains two singular points $q(2,0)$ and $q(2,\pi)$. At these points two simple singular curves $q(t^\pm(s),s)$ branches out. Here $t^\pm(s)=2\pm\sqrt{-2\sin s}$, $s\in [\pi, 2\pi]$.  
The simple singular curves touch each other at singular points, and actually they form a unique closed smooth curve. Moreover, at both singular points the genericity assumption is verified, so at these points we expect to have two germs of self-intersections.  

In order to complete the picture we apply the criterion of Proposition {\ref{Pr4}}. Assume that $a$ and $b$ are two distinct numbers on $[0,2\pi]$. Notice that  $v(s)=(\cos s, \sin s,0)$  for all $s$. Therefore the pairs $a$ and $b$, which may generate self-intersections, necessarily  satisfy  the relation $\cos a -\cos b=0$.  The non-trivial pairs of solutions of this equation are given by the pairs $(a,b)$ where $a=-b\;{\rm mod}\;2\pi$. Substituting this condition  into (\ref{self}) and performing all necessary simplifications we find the following non-trivial pairs of solutions:
$$
a=-b\; {\rm mod}\;2\pi,\qquad \tau_1=-2,\;\tau_2=2.
$$
Thus the surface in question intersects itself along the segment connecting the singular points $q(2,0)$ and $q(2,\pi)$. 
Moreover, for $a=\frac{\pi}{2}$ and $b=\frac{3\pi}{2}$ any pair of numbers $\tau_i$ satisfying $\tau_1-\tau_2+4=0$ is 
a solution of equation  (\ref{self}). Therefore the whole line passing throw the points $\gamma(\frac{\pi}{2})=(0,-2,0)$ and $\gamma(\frac{3\pi}{2})=(0,2,0)$ is a line of self-intersection.
In Fig.5 a) we show how the sub-Riemannian minimal surface described in this example looks like, and in Fig. 5 b) its singular set (fat curves).
}
\end{Ex}

\ppotRR{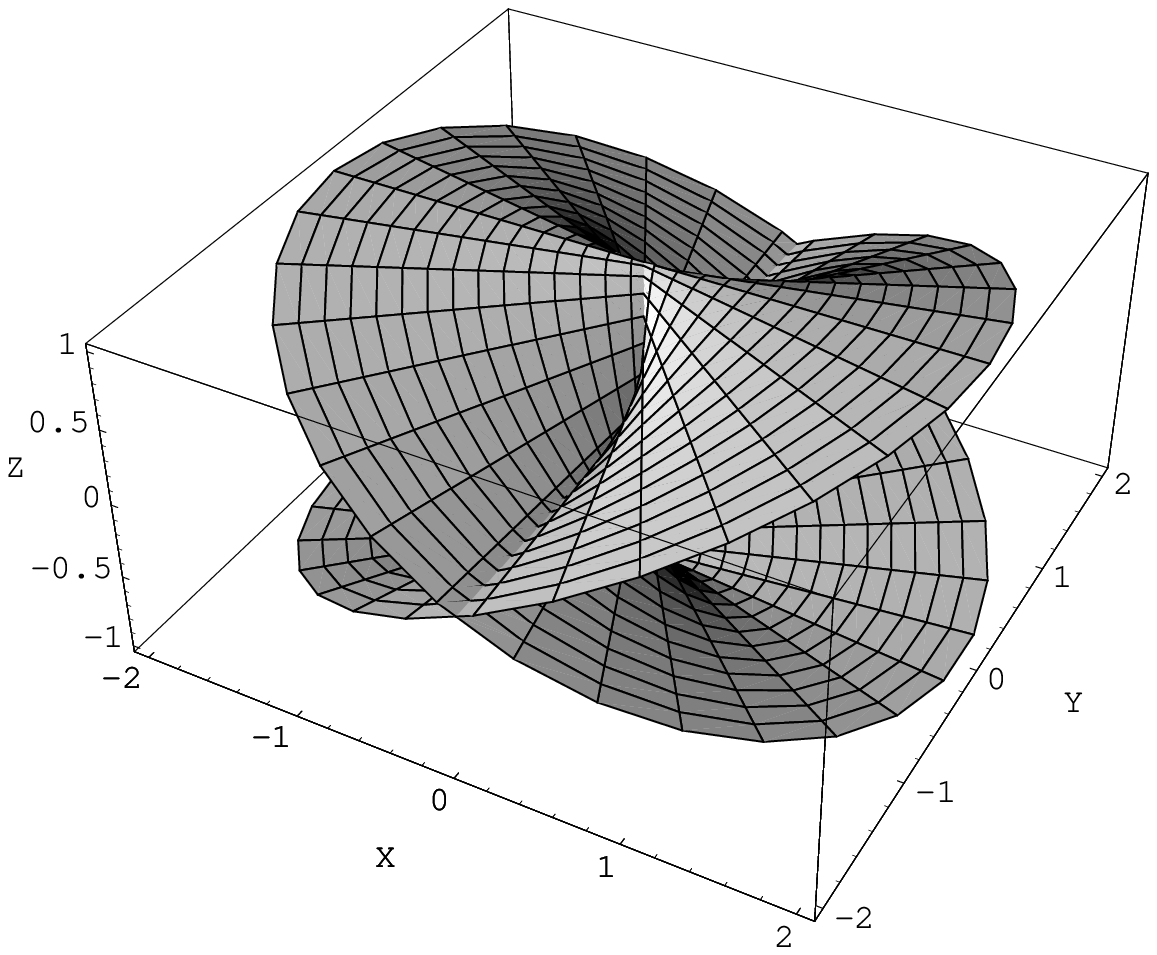}{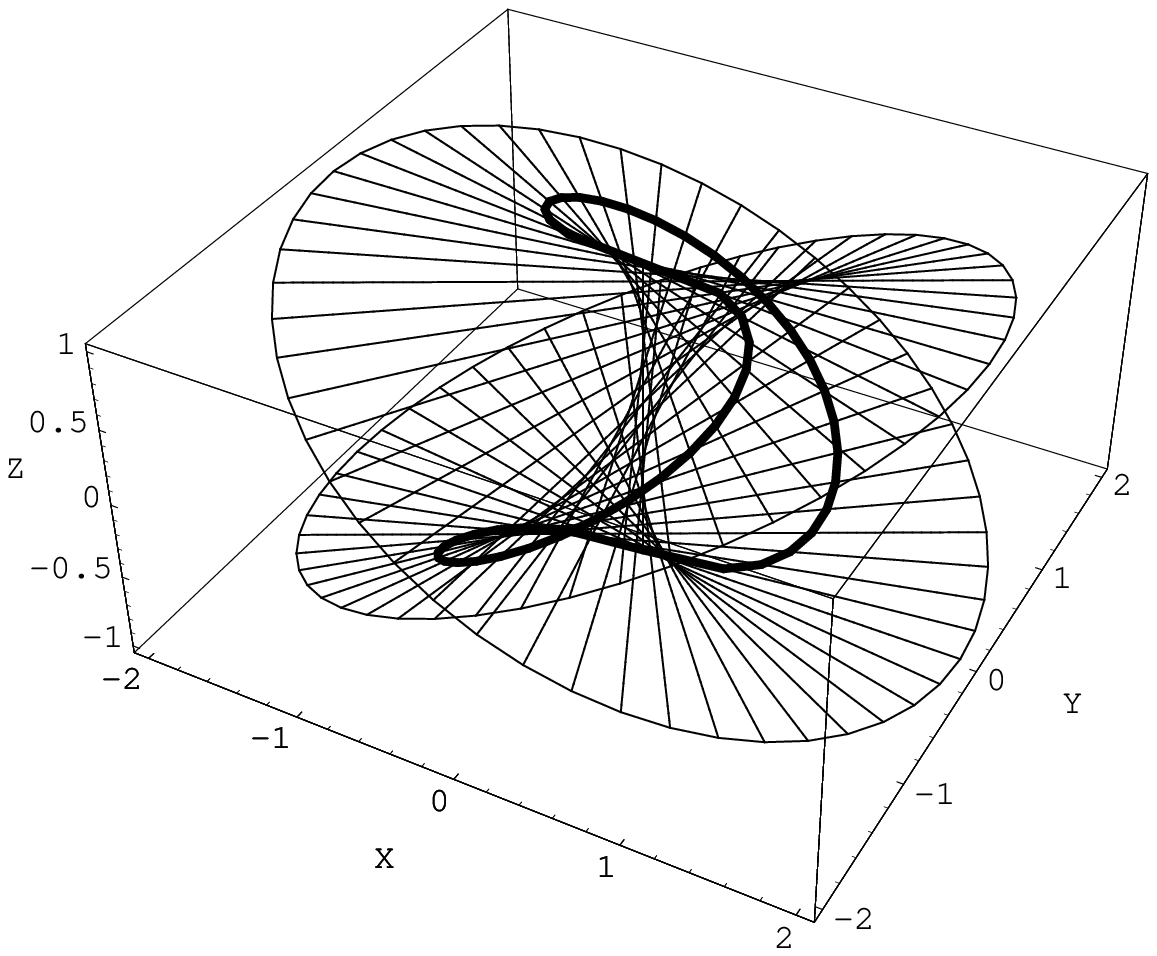}{Illustration to Example \ref{ex3}}{6.0}{6.0}

\end{document}